\documentclass[11pt, reqno]{amsart}

\usepackage{amssymb}
\usepackage{amsmath}

\usepackage{enumerate}

\usepackage{fullpage}
\def\stackunder#1#2{\mathrel{\mathop{#2}\limits_{#1}}}

\allowdisplaybreaks

\begin{document}

\newtheorem{theorem}{Theorem}[section]
\newtheorem{proposition}[theorem]{Proposition}
\newtheorem{definition}[theorem]{Definition}
\newtheorem{corollary}[theorem]{Corollary}
\newtheorem{lemma}[theorem]{Lemma}
\newtheorem{remark}[theorem]{Remark}
\newtheorem{example}[theorem]{Example}
\newtheorem{conjecture}[theorem]{Conjecture}

\title[Scaling limit of stochastic dynamics]{Scaling limit of
stochastic dynamics
in classical continuous systems}
\author{Martin Grothaus, Yuri G.~Kondratiev, Eugene Lytvynov,
Michael R\"ockner}

\address{Martin Grothaus, BiBoS and Mathematics Department,
Bielefeld University, 33615 Bielefeld, Germany and
IAM, Bonn University, 53115 Bonn, Germany. \newline
{\rm \texttt{Email: grothaus@wiener.iam.uni-bonn.de},
\texttt{URL: http://wiener.iam.uni-bonn.de/$\sim$grothaus/}}. \newline
Yuri G.~Kondratiev, BiBoS and Mathematics Department,
Bielefeld University, 33615 Bielefeld, Germany and
Inst.~Math., NASU, 252601 Kiev, Ukraine. \newline
{\rm \texttt{Email: kondrat@mathematik.uni-bielefeld.de}}. \newline
Eugene Lytvynov, BiBoS, Bielefeld University, 33615 Bielefeld, Germany
and IAM, Bonn University, 53115 Bonn, Germany. \newline
{\rm \texttt{Email: lytvynov@wiener.iam.uni-bonn.de}}. \newline
Michael R\"ockner, BiBoS and Mathematics Department,
Bielefeld University, 33615 Bielefeld, Germany. \newline
{\rm \texttt{Email: roeckner@mathematik.uni-bielefeld.de}
}}

\date{\today}

\thanks{We thank Sergio Albeverio, Hans Otto Georgii, Tobias Kuna, and Herbert
Spohn for discussions and helpful comments.
We would also like to thank a very conscientious referee for a careful
reading of the manuscript and making very useful comments and suggestions.
\newline
Financial support of the DFG through SFB 256 Bonn and the DFG-Forschergruppe
``Spectral Analysis, Asymptotic Distributions, and Stochastic Dynamics''
is gratefully acknowledged.}


\begin{abstract}
We investigate a scaling limit of gradient stochastic dynamics
associated to Gibbs states in classical continuous systems on
${\mathbb R}^d, \, d \ge 1$. The aim is to derive
macroscopic quantities from a given micro- or mesoscopic system.
The scaling we consider has been investigated
in \cite{Br80}, \cite{Ro81}, \cite{Sp86}, and \cite{GP86}, under
the assumption that the underlying potential is in $C^3_0$ and
positive. We prove that the Dirichlet
forms of the scaled stochastic dynamics converge on a core of
functions to the Dirichlet form of a
generalized Ornstein--Uhlenbeck process. The proof is based on the
analysis and geometry on the configuration space which was
developed in \cite{AKR98a}, \cite{AKR98b}, and works for general
Gibbs measures of Ruelle type. Hence, the underlying potential may
have a singularity at the origin, only has to be bounded from
below, and may not be compactly supported. Therefore, singular
interactions of physical interest are covered, as e.g.~the one given
by the Lennard--Jones
potential, which is studied in the theory of fluids. Furthermore,
using the Lyons--Zheng decomposition we give a simple proof for the
tightness of the scaled processes. We also prove that the
corresponding generators, however, do not converge in the $L^2$-sense.
This settles a conjecture formulated in
\cite{Br80}, \cite{Ro81}, \cite{Sp86}.
\end{abstract}

\maketitle

\section{Introduction}

The stochastic dynamics $({\bf X}(t))_{t \ge 0}$ of a classical
continuous system is an infinite dimensional diffusion process
having a Gibbs measure $\mu$, e.g.~of the type studied by Ruelle in
\cite{Rue69}, as an invariant measure. Physically, it
describes the stochastic dynamics of Brownian particles which are
interacting via the gradient of a pair-potential $\phi$. Since
each particle can move through each position in space, the system
is called continuous and is used for modelling gas and
fluid. For realistic models which can be described by
these stochastic dynamics, e.g.~suspensions, we refer to
\cite{Sp86}.

Since these dynamics are stochastic, they have to be interpreted
as mesoscopic processes. The aim of analyzing scaling limits, in general, is
to derive from micro- or mesoscopic systems macroscopic statements
and quantities. The type of scaling to study depends on which
features of a given system one is interested in, see
e.g.~\cite{Br80}, \cite{KL99}, \cite{Sp91}.

The scaling we consider in this paper has been
investigated in \cite{Br80} and \cite{Ro81}. In his Doctor-thesis,
\cite{Br80}, T.~Brox has given some heuristic arguments for
non-convergence in law of the scaled process and has conjectured that
there is no limiting Markov process. However, assuming the convergence
of the generators of the scaled stochastic dynamics averaged over
time, cf.~Conjecture \ref{le2000} below, H.~Rost
has given some heuristic arguments in \cite{Ro81} for the existence of a limiting
generalized Ornstein--Uhlenbeck process, which, of course,
contradicts the statement of Brox.
A fundamental and celebrated paper on this problem is due to
H.~Spohn \cite{Sp86}. Assuming that the underlying
potential is smooth, compactly supported and positive,
there the author describes a proof of Conjecture \ref{le2000}
within the proof of his main theorem (see, however, the remark on page
4 of \cite{Sp86}, and Proposition 2 therein, concerning the restriction
$d \le 3$). Another
approach has been proposed in \cite{GP86}. The idea of M.~Z.~Guo
and G.~Papanicolaou has been to prove convergence of the
corresponding resolvent. As remarked by themselves, at that time the
authors did not have an appropriate infinite dimensional analysis
and geometry at their disposal, and therefore their considerations
have been on a non-rigorous level.

After these contributions, for a long time there has been no progress in
this problem. Recently, however, some new techniques have been introduced.
In \cite{AKR98a}, \cite{AKR98b} an infinite
dimensional analysis and geometry on the configuration space was
developed. In this paper we shall make use of these concepts in order
to tackle the problem described above again.

The stochastic dynamics $({\bf X}(t))_{t \ge 0}$ of a classical
continuous system takes values in the configuration space
\begin{eqnarray*}
\Gamma := \{ \gamma \subset {\mathbb R}^{d} | \, |\gamma \cap K| <
\infty \,\,\, \mbox{for any compact} \, K \subset {\mathbb R}^{d} \},
\end{eqnarray*}
and informally solves the following infinite system of stochastic
differential equations:
\begin{align}\label{eq0621}
dx(t) & = -\beta \sum_{\stackunder{y(t) \neq x(t)} {y(t)
\in {\bf X}(t)}} \nabla \phi(x(t) - y(t)) \,dt + \sqrt{2} \,dB^x(t),  &
x(t) \in {\bf X}(t), \nonumber \\ {\bf X}(0) & = \gamma,  &\gamma
\in \Gamma,
\end{align}
where $(B^x)_{x \in \gamma}$ is a sequence of
independent Brownian motions. The study of such diffusions has
been initiated by R.~Lang \cite{La77} (see also \cite{Sh79}),
who considered the case
$\phi \in C_0^3({\mathbb R}^d)$ using finite dimensional
approximations and stochastic
differential equations. More singular $\phi$, which are of
particular interest in Physics, as e.g.~the Lennard--Jones
potential, have been treated by H.~Osada, \cite{Os96}, and M.~Yoshida,
\cite{Y96} (see also \cite{Ta97}, \cite{FRT00} for
the hard core case). Osada and Yoshida were the first to
use Dirichlet forms for the construction of such processes.
However, they could not write down the corresponding
generators or martingale problems explicitly, hence could not prove
that their processes actually solve (\ref{eq0621}) weakly.
This, however, was proved in \cite{AKR98b} by showing an integration by
parts formula for the respective Gibbs measures. Thus the latter work
became the starting point of this paper. In \cite{AKR98b},
also Dirichlet forms were used and all
constructions were designed to work particularly for singular potentials
of the above mentioned type, see Theorem \ref{th57} below.
Additionally, and this is essential for our considerations, an explicit
expression for the corresponding generator and martingale
problem was provided, which shows that the process in \cite{AKR98b}
indeed solves (\ref{eq0621}) in the weak sense.

The scaled process $({\bf
X}_{\epsilon}(t))_{t \ge 0}$ studied in this paper is defined by
\begin{eqnarray*}
{\bf X}_{\epsilon}(t) := S_{out, \epsilon} (S_{in, \epsilon} ({\bf X}
(\epsilon^{-2}t))), \qquad t \ge 0, \quad \epsilon > 0,
\end{eqnarray*}
and we are interested in the scaling limit for $\epsilon \to 0$.
The first scaling $S_{in, \epsilon}$ scales the position of the
particles inside the configuration space as follows:
\begin{eqnarray*}
\Gamma \ni \gamma \mapsto S_{in, \epsilon}(\gamma)
:= \{ \epsilon \, x | \, x \in \gamma \} \in \Gamma, \qquad \epsilon > 0.
\end{eqnarray*}
Hence, for small $\epsilon > 0$ this scaling concentrates the particles
towards the origin. The second scaling $S_{out, \epsilon}$ leads us out
of the configuration space and is given by
\begin{eqnarray*}
\Gamma \ni \gamma \mapsto S_{out, \epsilon}(\gamma) :=
\epsilon^{d/2}\Big(\gamma - \rho^{(1)}_{\tilde{\mu}_\epsilon} \,dx
\Big) \in {\mathcal D}',
\end{eqnarray*}
where ${\mathcal D}'$ is the dual space of ${\mathcal D}
:= C_0^{\infty}({\mathbb R}^d)$. In the second
scaling we first center the configuration $\gamma$ by subtracting
the first correlation measure $\rho^{(1)}_{\tilde{\mu}_\epsilon}
\,dx$ of the Gibbs measure $\tilde{\mu}_\epsilon := S_{in,
\epsilon}^{\ast}\mu$. Furthermore, we scale the mass of the
particles by $\epsilon^{d/2}$ to avoid divergence of the total
mass at the origin as $\epsilon \to 0$.

We start with constructing the Dirichlet form ${\mathcal
E}_\epsilon$, the generator $H_\epsilon$ and the semi-group $(T_{\epsilon,
t})_{t \ge 0}$ associated to $({\bf X}_{\epsilon}(t))_{t \ge 0}$.
These objects are images of the Dirichlet form,
generator, and semi-group, respectively, which are
associated to the original stochastic
dynamics $({\bf X}(t))_{t \ge 0}$, see Theorem \ref{th58} below.

The first convergence we show is the following, see Theorem
\ref{th11}. We prove that
\begin{eqnarray}\label{eq1111}
\lim_{\epsilon \to 0} {\mathcal E}_{\epsilon}(F,G) = {\mathcal
E}_{\nu_{\mu}}(F,G),
\end{eqnarray}
for all smooth cylinder functions $F, G \in {\mathcal F}C_b^{\infty}({\mathcal
D}, {\mathcal D}')$. The limit Dirichlet form ${\mathcal E}_{\nu_{\mu}}$
is defined on $L^2({\mathcal D}', \nu_\mu)$ with $\nu_\mu$ being white noise,
and associated to a generalized
Ornstein--Uhlenbeck process $({\bf X}(t))_{t \ge 0}$ solving the
stochastic differential equation
\begin{eqnarray}\label{eq0121}
d{\bf X}(t,x) =
\frac{\rho_{\phi}^{(1)}(\beta,1)}{\chi_{\phi}(\beta) } \Delta {\bf
X}(t,x) \,dt + \sqrt{2 \, \rho_{\phi}^{(1)}(\beta,1)} \,d{\bf W}(t,x),
\end{eqnarray}
where $({\bf W}(t))_{t \ge 0}$ is a Brownian motion in ${\mathcal
D}^{\prime}$ with covariance operator $- \Delta$. The coefficient
$\rho_{\phi}^{(1)}(\beta,1) / \chi_{\phi}(\beta) $ is called the bulk
diffusion coefficient and $\beta$ is the inverse temperature. The
convergence (\ref{eq1111}) determines the limit process uniquely,
see Remark \ref{rm1111}(i), and requires only very weak
assumptions. The interaction potential $\phi$ only has to be
stable (S) and we have to assume the LA-HT, i.e.~the low activity
high temperature regime (see below for precise definitions).
A basic ingredient in the proof is the
convergence of the image measures $\mu_\epsilon := S_{out,
\epsilon}^{\ast} S_{in, \epsilon}^{\ast} \mu$ to the Gaussian
white noise measure $\nu_\mu$ as $\epsilon \to 0$, see Theorem
\ref{th9}. The latter fact has been proved by T.~Brox, \cite{Br80}.

The convergence in terms of the Dirichlet forms, however,
up to this point has no
probabilistic interpretation. Hence, we also study convergence in
law of the scaled processes. By ${\bf P}^{\epsilon}$ we
denote the law of the scaled equilibrium processes, i.e., the law
of the scaled process starting with a distribution equal to
the equilibrium measure $\mu_\epsilon$. Then, in Theorem
\ref{pr777} we prove that the family $({\bf
P}^{\epsilon})_{\epsilon > 0}$ is tight. This has been shown before by
T.~Brox, \cite{Br80}, and H.~Spohn, \cite{Sp86}, for smooth compactly
supported potentials. Our proof, again, works under quite weak
assumptions on the potential. We only need
conditions which ensure the existence of the original stochastic
process and have to assume the LA-HT regime. In the proof we
use the well-known Lyons-Zheng decomposition,
\cite{LZ88}, \cite{LZ94}, of
the scaled process and the Burkholder--Davies--Gundy
inequalities in order to establish the required estimate of the
increments. Since the state space of the scaled process is a space
of distributions, we first prove tightness in a weak
sense. Then, via some Hilbert--Schmidt embeddings, we find a negative,
weighted Sobolev spaces ${\mathcal H}_{-m}$ as state space
such that the family $({\bf P}^{\epsilon})_{\epsilon > 0}$
is tight on $C([0, \infty), {\mathcal H}_{-m})$.

It remains to prove that all accumulation points coincide with
the generalized Ornstein--Uhlenbeck
process $({\bf X}(t))_{t \ge 0}$ above. A well-known method to
identify the limit is based on considering the associated martingale problem.
More precisely, if we could prove that all accumulation points of $({\bf
P}^{\epsilon})_{\epsilon > 0}$ satisfy the martingale
problem for the generator $H$ associated to equation (\ref{eq0121})
with initial condition $\nu_\mu$, then a (slight
modification of a) uniqueness result of
R.~Holley and D.~Stroock \cite{HS78} implies that all these
accumulation points coincide.

The obvious first idea to prove that all limit points solve the
martingale problem for $H$ is to try to prove strong convergence of
$H_\epsilon \to H$ as $\epsilon \to 0$.
In \cite{Br80}, \cite{Ro81}, and \cite{Sp86}
it has, however, been conjectured that, in general, the difference
\begin{eqnarray*}
\parallel (H - H_\epsilon)F \parallel_{L^2(\mu_\epsilon)},
\qquad F \in {\mathcal F}C_b^{\infty}({\mathcal D}, {\mathcal D}'),
\end{eqnarray*}
does not tend to zero as $\epsilon \to 0$. In Theorem \ref{th2001}
we prove that this conjecture is indeed true.
The proof is quite an elaborate task and is done via a (mathematically
rigorous) high temperature expansion.
A basic tool for this is provided by Theorem
\ref{pr200}, where we derive explicit formulas for the derivative
of the correlation functions with respect to the inverse
temperature $\beta$ using the so-called $K$-transform from \cite{KK99a},
and by Theorem \ref{dfdrr}, where we prove a
coercivity identity for Gibbs measures.

It turns out that for the above described identification of the accumulation points of
$({\bf P}^{\epsilon})_{\epsilon > 0}$, however, a weaker
convergence of the generators is sufficient. In Theorem \ref{th2000}
we prove convergence in law under the assumption that Conjecture
\ref{le2000} is true, i.e., under the assumption that the
generators converge in time average.

To complete the program also from a purely probabilistic point of
view, it remains to prove Conjecture \ref{le2000} in physically
relevant models. This will be the subject of future work.

The progress achieved in this paper may be summarized
by the following core results:
\begin{itemize}
\item Convergence of Dirichlet forms is shown, see Remark \ref{rm0002}.

\item The tightness result as in \cite{Br80}, \cite{Sp86} is generalized, see
Remark \ref{rm1111}.

\item Conjecture on non-convergence of generators is proved.

\item A mathematically rigorous high temperature expansion of
all correlation functionals is developed (up to second order in $\beta = 1/T$).

\item All above results apply to physically relevant potentials, in particular
singularities at the origin, non-trivial negative part, and infinite range are allowed.
\end{itemize}

Hypotheses on the potential are weakened not for the sake of generality, but in order
to cover the physically relevant potentials (as e.g.~ Lennard--Jones
potential).

\section{Gibbs states of classical continuous systems}

\subsection{Configuration space and Poisson measure}

Let ${\mathbb R}^{d}, d \ge 1$, be equipped with the norm
$|\cdot|_{{\mathbb R}^{d}}$
given by the Euclidean scalar product $(\cdot, \cdot)_{{\mathbb
R}^{d}}$. By ${\mathcal B}({\mathbb R}^{d})$ we denote the
corresponding Borel $\sigma$-algebra. ${\mathcal O}_c({\mathbb
R}^{d})$ denotes the system of all open sets in ${\mathbb
R}^{d}$, which have compact closure. The Lebesgue measure on the
measurable space $({\mathbb R}^{d}, {\mathcal B}({\mathbb R}^{d}))$ we
denote by $dx$.

The {\it configuration space} $\Gamma$ over
${\mathbb R}^{d}$ is defined by
\begin{eqnarray*}
\Gamma := \{ \gamma \subset {\mathbb R}^{d} | \, |\gamma \cap K| < \infty \,
\,\, \mbox{for any compact} \, K \subset {\mathbb R}^{d} \}.
\end{eqnarray*}
Here $|A|$ denotes the cardinality of a set $A$.
Via the identification of $\gamma \in \Gamma$ with
$\sum_{x \in \gamma} \varepsilon_{x} \in {\mathcal M}_p({\mathbb R}^{d})$,
where $\varepsilon_{x}$ denotes the Dirac measure in $x \in
{\mathbb R}^{d}$, $\Gamma$ can be considered as a subset of the
set ${\mathcal M}_p({\mathbb R}^{d})$ of all positive Radon measures
on ${\mathbb R}^{d}$. Hence $\Gamma$ can be topologized by the
vague topology, i.e., the topology generated by maps
\begin{eqnarray*}
\gamma \mapsto \, \langle f,\gamma  \rangle \, := \int_{{\mathbb
R}^{d}} f(x) \,d\gamma(x) = \sum_{x \in \gamma} f(x)
\end{eqnarray*}
where $f \in C_{0}({\mathbb R}^{d})$,
the set of continuous functions on ${\mathbb R}^{d}$ with compact support.
We denote by ${\mathcal B}({\Gamma})$ the corresponding Borel $\sigma$-algebra.

For a given $z>0$ (activity  parameter), let $ \pi_{z} $ denote
the Poisson measure on $(\Gamma, {\mathcal B}(\Gamma))$ with intensity
measure $z\,dx$. This measure is characterized via its
Fourier transform
\begin{eqnarray*}
\int_{\Gamma} \exp(i\langle f,\gamma  \rangle) \,d\pi_{z}(\gamma)
= \exp\Big{(} z \int_{{\mathbb R}^{d}} (\exp(if(x)) - 1) \,dx \Big{)},
\qquad f \in {\mathcal D},
\end{eqnarray*}
where ${\mathcal D} := C^{\infty}_0({\mathbb R}^d)$, the set of smooth functions
on ${\mathbb R}^d$ with compact support.

\subsection{Gibbs measures in the LA-HT regime}

Let $\phi$ be a symmetric pair potential, i.e., a measurable function
$\phi: {\mathbb R}^d \to {\mathbb R} \cup \{ \infty \}$ such that
$\phi(x) = \phi(-x)$. For $\Lambda \in {\mathcal O}_c({\mathbb R}^{d})$ the
conditional energy $E^{\phi}_{\Lambda}: \Gamma  \to {\mathbb R} \cup \{ \infty \}$
with empty boundary condition is defined by
\begin{eqnarray*}
E_{\Lambda}^{\phi}(\gamma) :=
\sum\limits_{ \{x,y\} \subset \gamma_{\Lambda}} \phi(x-y)
= E_{\Lambda}^{\phi}(\gamma_\Lambda)
\end{eqnarray*}
where $\gamma_{\Lambda} := \gamma \cap \Lambda$ and the sum over
the empty set is defined to be zero.

For every $r = (r_1, \ldots, r_d) \in {\mathbb Z}^d$ we define a cube
\begin{eqnarray*}
Q_r = \Big{\{} x \in {\mathbb R}^d \, \Big{|} \, r_i - 1/2 \le x_i < r_i + 1/2 \Big{\}}.
\end{eqnarray*}
These cubes form a partition of ${\mathbb R}^d$. For any $\gamma \in \Gamma$ we set
$\gamma_r := \gamma_{Q_r}, \, r \in {\mathbb Z}^d$. Additionally, we introduce for
$n \in {\mathbb N}$ a cube $\Lambda_n$ with side length $2n -1$ centered at the origin
in ${\mathbb R}^d$.

Let us recall some standard assumptions from Statistical Mechanics. For
our results we have to require some of the following conditions.

\begin{description}
\item[(SS)] ({\it superstability}) There exist $A(\phi) > 0, \, B(\phi) \ge 0$ such that, if $\gamma =
\gamma_{\Lambda_n}$ for some $n \in {\mathbb N}$, then
\begin{eqnarray*}
E_{\Lambda_n}^{\phi}(\gamma) \, \ge \, \sum_{r \in {\mathbb Z}^d}
\Big{(} A(\phi) |\gamma_r|^2 - B(\phi) |\gamma_r| \Big{)}.
\end{eqnarray*}
\end{description}
(SS) obviously implies:
\begin{description}
\item[(S)] ({\it stability})
For any $\Lambda \in {\mathcal O}_c({\mathbb R}^{d})$ and for all
$\gamma \in \Gamma$ we have
\begin{eqnarray*}
E_{\Lambda}^{\phi}(\gamma) \, \ge \, -B(\phi) |\gamma_{\Lambda}|.
\end{eqnarray*}
\end{description}
A consequence of (S), in turn, is, of course, that $\phi$ is bounded from below.
For $\beta \ge 0, z > 0$, let us define
\begin{eqnarray*}
C(\beta \phi, z) :=  \exp(2\beta B(\phi))
\int_{{\mathbb R}^{d}} | \exp(-  \beta \phi(x)) - 1 | \,
z \,dx.
\end{eqnarray*}
We also need
\begin{description}
\item[(UI)] ({\it uniform integrability}) We have:
\begin{eqnarray*}
C(\beta \phi, z) < \exp(-1).
\end{eqnarray*}
\end{description}
For a given potential $\phi$ the set of pairs $(\beta, z)$ such that
Condition (UI) holds is called {\bf LA-HT (low activity high temperature) regime},
see \cite{Rue63} and \cite{Mi67}. (UI) is stronger than (I)
({\it integrability}), i.e., $C(\beta \phi, z) < \infty$,
which is also called {\it regularity}, see e.g.~\cite{Rue69}.
\begin{description}
\item [(LR)] ({\it lower regularity}) There exists a decreasing positive function $a: {\mathbb N}
\to {\mathbb R}_+$ such that
\begin{eqnarray*}
\sum_{r \in {\mathbb Z}^d} a(\parallel r \parallel) < \infty
\end{eqnarray*}
and for any $\Lambda^{\prime}, \Lambda^{\prime \prime}$ which are finite unions
of cubes of the form $Q_r$ and disjoint,
\begin{eqnarray*}
W^{\phi}(\gamma^{\prime} \mid \gamma^{\prime \prime}) \ge
- \sum_{r^{\prime}, r^{\prime \prime} \in {\mathbb Z}^d} a(\parallel
r^{\prime} -  r^{\prime \prime} \parallel) \,
|\gamma^{\prime}_{r^\prime}| \,
|\gamma^{\prime \prime}_{r^{\prime \prime}}|,
\end{eqnarray*}
provided $\gamma^{\prime}
= \gamma^{\prime}_{\Lambda^{\prime}}, \, \gamma^{\prime \prime}
= \gamma^{\prime \prime}_{\Lambda^{\prime \prime}}$.
Here
\begin{eqnarray*}
W^{\phi}(\gamma' \mid \gamma'') :=
\sum\limits_{x \in \gamma', y \in \gamma''} \phi(x-y)
\end{eqnarray*}
is the interaction energy and $\parallel \cdot \parallel$ denotes
the maximum norm on ${\mathbb R}^d$.
\end{description}

On $(\Gamma, {\mathcal B}(\Gamma))$ we consider the finite volume Gibbs measures
$\mu_{\Lambda}$ in
$\Lambda \in {\mathcal O}_c({\mathbb R}^{d})$ with empty boundary condition:
\begin{eqnarray*}
d\mu_{\Lambda}(\gamma) := \frac{1}{Z_{\Lambda}} \exp (- \beta
E_{\Lambda}^{\phi} (\gamma_{\Lambda})) d\pi_{z}(\gamma),
\end{eqnarray*}
where $\beta \ge 0$ is the inverse temperature and
\begin{eqnarray*}
Z_{\Lambda} = \int_{\Gamma} \exp (- \beta E_{\Lambda}^{\phi}
(\gamma_{\Lambda})) d\pi_{z}(\gamma)
\end{eqnarray*}
is the partition function. Using (S) one easily proves that
it is finite. In e.g.~\cite{Mi67} and \cite{MM91} it has been proved that
in the LA-HT regime the weak limit
\begin{eqnarray}\label{eq2}
\lim_{{\Lambda \nearrow {\mathbb R}^{d}}} \mu_{\Lambda} = \mu
\end{eqnarray}
exists. Furthermore, it can be shown that $\mu$ is a Gibbs measure,
see \cite{Rue70} and \cite{K99}. The measure $\mu$ in
\eqref{eq2} we call Gibbs measure corresponding to $(\phi, \beta, z)$ and
the construction with empty boundary condition.

\subsection{$K$-transform and correlation functions}\label{ss23}

Next, we recall the definition of correlation functions using
the concept of the so-called $K$-transform, see e.g. \cite{KK99a},
\cite{Len73}, \cite{Len75a}, \cite{Len75b}.

Denote by $\Gamma_0$ the space of finite configurations over ${\mathbb R}^d$:
\begin{eqnarray*}
\Gamma_0 := \bigsqcup_{n=0}^\infty \Gamma_0^{(n)},\quad
\Gamma^{(0)}_0 := \{\varnothing\},\quad \Gamma_0^{(n)}:=\{\eta
\subset {\mathbb R}^d \mid |\eta|= n \}, \quad n \in {\mathbb N}.
\end{eqnarray*}
Let $\widetilde{{\mathbb R}^{d \times n}}= \big\{\,(x_1, \dots, x_n)
\in {\mathbb R}^{d  \times n} \mid x_i \ne x_j \quad \mbox{for}
\quad i \ne j \, \big\}$ and let $S^n$ denote the group of all permutations of
$\{1,\dots,n\}$. Through the natural bijection
\begin{eqnarray}\label{eq366}
\widetilde{{\mathbb R}^{d \times n}}/S^n \longleftrightarrow
\Gamma_0^{(n)}
\end{eqnarray}
one defines a topology on $\Gamma^{(n)}_0$. The space $\Gamma_0$ is
equipped then with the topology of disjoint union. Let ${\mathcal B}(\Gamma_0)$
denote the Borel $\sigma$-algebra on $\Gamma_0$.

A ${\mathcal B}(\Gamma_0)$-measurable function $G \colon \Gamma_0 \to
{\mathbb R}$ is said to have bounded support if there exist
$\Lambda \in {\mathcal O}_c({\mathbb R}^{d})$ and $N \in {\mathbb
N}$ such that $\mbox{supp}(G) \subset \bigsqcup_{n=0}^N
\Gamma_{0, \Lambda}^{(n)}$, where $\Gamma_{0, \Lambda}^{(n)} = \{\eta
\subset \Lambda \mid |\eta|= n \}$.

For any $\gamma \in {\Gamma}$ let $\sum_{\eta \Subset \gamma}$
denote the summation over all $\eta \subset \gamma$ such that
$|\eta| < \infty$. For a function
$G: \Gamma_0 \to {\mathbb R}$ , the $K$-transform of $G$ is
defined by
\begin{eqnarray}\label{eq9}
(KG)(\gamma):= \sum_{\eta \Subset \gamma} G(\eta)
\end{eqnarray}
for each $\gamma \in \Gamma$ such that at least one of the series
$\sum_{\eta \Subset \gamma} G^+(\eta)$ or $\sum_{\eta \Subset
\gamma} G^-(\eta)$ converges, where $G^{+} := \max \{ 0, G\}$ and
$G^{-} := -\min \{ 0, G\}$.

Let $\mu$ be a probability measure on $(\Gamma,{\mathcal B}(\Gamma))$. The
correlation measure corresponding to $\mu$ is defined by
\begin{eqnarray*}
\rho_\mu(A) := \int_{\Gamma}(K1_A)(\gamma) \,d\mu(\gamma), \qquad A
\in {\mathcal B}(\Gamma_0).
\end{eqnarray*}
$\rho_\mu$ is a measure on $(\Gamma_0, {\mathcal B} (\Gamma_0))$
(see \cite{KK99a} for details, in particular, measurability issues).

Let $G \in L^1(\Gamma_0, {\mathcal B}(\Gamma_0), \rho_\mu)$,
then $ \parallel KG \parallel_{L^1(\mu)} \le \parallel K|G| \parallel_{L^1(\mu)}
= \parallel G \parallel_{L^1(\rho_\mu)}$,
hence $KG \in L^1(\Gamma, {\mathcal B}(\Gamma), \mu)$ and $KG(\gamma)$
is for $\mu$-a.e.~$\gamma \in \Gamma$ absolutely convergent.
Moreover, then obviously
\begin{eqnarray}\label{eq302}
\int_{\Gamma_0} G(\eta) \,d\rho_\mu(\eta) =
\int_{\Gamma}(KG)(\gamma) \,d\mu(\gamma),
\end{eqnarray}
see \cite{KK99a}, \cite{Len75a}, \cite{Len75b}.

The Lebesgue--Poisson measure $\lambda$ on $(\Gamma_0,{\mathcal B}(\Gamma_0))$
with activity parameter $z > 0$ is defined by
\begin{eqnarray*}
\lambda_z := \delta_\emptyset+\sum_{n=1}^\infty \frac{z^n}{n!}\,dx^{\otimes n},
\end{eqnarray*}
where $dx^{\otimes n}$ is defined via the bijection (\ref{eq366}).

For the Gibbs measure $\mu$ in the LA-HT regime corresponding
to $\phi$ satisfying (S) and the construction with empty boundary
condition, the correlation
measure $\rho_\mu$ is absolutely continuous with respect to the
Lebesgue-Poisson measure, see e.g.~\cite{Rue63} and \cite{Mi67}.
Its Radon-Nikodym derivative
\begin{eqnarray*}
\rho_{\mu}(\eta) := \frac{d\rho_\mu}{d\lambda_z}(\eta), \qquad \eta \in
\Gamma_0,
\end{eqnarray*}
w.r.t.~$\lambda_z$ we denote by the same symbol and the functions
\begin{eqnarray}\label{eq367}
\rho_{\mu}^{(n)}(x_1,\dots,x_n) := \rho_{\mu}(\{x_1,\dots,x_n\}),
\quad x_1,\dots,x_n \in {\mathbb R}^d, \,\, x_i \neq x_j
\,\, \mbox {if} \,\, i \neq j,
\end{eqnarray}
are called the $n$-th order correlation functions of the measure
$\mu$. Furthermore, the correlation functions can be expressed as
functions of the underlying potential $\phi$, inverse temperature
$\beta$ and activity $z$, i.e.,
$\rho_\mu = \rho_\phi(\beta, z)$, see e.g.~\cite{Rue63}, \cite{Mi67}.
Hence, due to the translation invariance of the pair
interaction, the correlation functions as well as the
Gibbs measure $\mu$ are translation invariant. In particular,
$\rho_{\phi}^{(1)}(\beta, z)$ does not depend on $x_1 \in {\mathbb
R}^d$.

Additionally, for these functions the so-called Ruelle bound
holds: for fixed $\beta \ge 0, z > 0,$
there exists a constant $\xi > 0$ such that for
all $n$ and $x_1,\dots,x_n \in {\mathbb R}^d, \,
x_i \neq x_j$ for $i \neq j$, we have
\begin{eqnarray}\label{eq18}
\rho_\phi^{(n)}(\beta, z, x_1,\dots,x_n) \le \xi^n,
\end{eqnarray}
see \cite{Rue63}.
Using this bound one, in particular, gets that all local moments
of $\mu$ are finite:
\begin{eqnarray}\label{eq19}
\int_\Gamma |\gamma_\Lambda|^n \,d\mu(\gamma) < \infty \qquad
\forall n \in {\mathbb N}, \qquad \Lambda \in {\mathcal O}_c({\mathbb R}^{d}).
\end{eqnarray}

\section{Dirichlet forms, their generators, and corresponding
stochastic dynamics}\label{s44}

Here we recall the analysis and geometry on configuration space developed
in \cite{AKR98a} and \cite{AKR98b}.

Let $T_x({\mathbb R}^d)={\mathbb R}^d$ denote the tangent space to
${\mathbb R}^d$ at a point $x\in {\mathbb R}^d$. The tangent space
to $\Gamma$ at a point $\gamma\in\Gamma$ is defined as the Hilbert
space
\begin{eqnarray*}
T_\gamma(\Gamma) := L^2({\mathbb R}^d\to T{\mathbb
R}^d,\gamma )=\bigoplus_{x\in\gamma}T_x({\mathbb R}^d).
\end{eqnarray*}
Thus, each $V(\gamma )\in T_\gamma (\Gamma )$ has the form
$V(\gamma )=(V(\gamma,x ))_{x\in \gamma }$, where $V(\gamma,x )\in
T_x({\mathbb R}^d)$, and
\begin{eqnarray*}
\| V(\gamma )\| _{T_\gamma(\Gamma)}
^2=\sum_{x\in \gamma }\|V(\gamma,x )\|_{T_x({\mathbb
R}^d)}^2=\sum_{x\in \gamma }\|V(\gamma,x )\|_{{\mathbb R}^d}^2.
\end{eqnarray*}

Let $\gamma \in \Gamma $ and $x\in \gamma $. We denote by
${\mathcal O}_{\gamma ,x}$  an arbitrary open neighborhood of $x$
in $X$ such that ${\mathcal O}_{\gamma,x}\cap(\gamma \setminus
\{x\})=\varnothing$.  Now, for a function $F\colon\Gamma \to
{\mathbb R}$, $\gamma \in \Gamma $, and $x\in\gamma $, we define a
function $F_{x}(\gamma ,\cdot )\colon{\mathcal O}_{\gamma ,x}\to
{\mathbb R}$ by
\begin{eqnarray*}{\mathcal O}_{\gamma,x}\ni y\mapsto
F_x(\gamma,y) := F(\gamma-\varepsilon_x+\varepsilon_y)\in{\mathbb
R}.
\end{eqnarray*}

We say that a function $F\colon \Gamma \to {\mathbb R}$ is
differentiable at $\gamma \in \Gamma $ if, for each $x\in \gamma
$, the function $F_x(\gamma ,\cdot )$ is differentiable at $x$ and
\begin{eqnarray*}
\nabla ^\Gamma F(\gamma ) := (\nabla_x  F(\gamma ))_{x\in \gamma
}\in T_\gamma( \Gamma ),
\end{eqnarray*}
where
\begin{eqnarray*}
\nabla_x  F(\gamma ) := \nabla _yF_x(\gamma ,y){\big |}_{y=x}.
\end{eqnarray*}
Evidently, this definition is independent
of the choice of the set ${\mathcal O}_{\gamma,x}$. We call
$\nabla^\Gamma F(\gamma)$ the gradient of $F$ at $\gamma \in \Gamma$.

We define a set of smooth cylinder functions
${\mathcal F}C_b^{\infty}({\mathcal D}, \Gamma)$ as the set of all functions on
$\Gamma$ of the form
\begin{eqnarray}\label{eq54}
\gamma \mapsto F(\gamma) = g_F(\langle f_1, \gamma \rangle, \ldots,
\langle f_N, \gamma \rangle),
\end{eqnarray}
where $f_1, \ldots, f_N \in {\mathcal D}$ and $g_F \in
C^{\infty}_b({\mathbb R}^N)$. Clearly, ${\mathcal F}C_b^{\infty}({\mathcal
D}, \Gamma)$ is dense in $L^p(\mu), p \ge 1$. Any function $F$
of the form (\ref{eq54}) is differentiable at each point
$\gamma\in\Gamma$, and its gradient is given by
\begin{eqnarray}\label{eq78}
(\nabla^\Gamma F)(\gamma, x) = \sum_{j=1}^{N} \partial_j g_F
(\langle f_1, \gamma \rangle, \ldots,
\langle f_N, \gamma \rangle) \nabla f_j(x), \quad
\gamma \in \Gamma, \, x \in \gamma,
\end{eqnarray}
where $\partial_j$ denotes the partial derivative w.r.t.~the
$j$-th variable.
For $F, G \in {\mathcal F}C_b^{\infty}({\mathcal D}, \Gamma)$ we define
\begin{eqnarray*}
{\mathcal E}_{\mu}^{\Gamma}(F,G) := \int_{\Gamma} (\nabla^{\Gamma}F(\gamma),
\nabla^{\Gamma}G(\gamma))_{T_{\gamma}(\Gamma)} \,d \mu(\gamma).
\end{eqnarray*}
Gibbs measures $\mu$ in the LA-HT regime corresponding to stable potentials
and the construction with empty boundary condition have all local moments finite,
see (\ref{eq19}). Thus, for such measures with the help of
(\ref{eq78}) we have $(\nabla^{\Gamma}F(\gamma),
\nabla^{\Gamma}G(\gamma))_{T_{\gamma}(\Gamma)} \in L^1(\mu)$.
Furthermore, the gradient respects $\mu$-classes ${\mathcal
F}C_b^{\infty}({\mathcal D}, \Gamma)^{\mu}$ determined by ${\mathcal
F}C_b^{\infty}({\mathcal D}, \Gamma)$, see e.g.~\cite{Ro98}, \cite{MaRo00}.
Hence, $({\mathcal
E}_{\mu}^{\Gamma}, {\mathcal F}C_b^{\infty}({\mathcal D}, \Gamma))$ is a
densely defined, positive definite, symmetric bilinear form on
$L^2(\mu)$.

To ensure closability of this bilinear form we have
to assume further properties of the potential $\phi$ :

\begin{description}
\item[(D)] ({\it differentiability})
The function $\exp(-\phi)$ is weakly differentiable on ${\mathbb
R}^d$, $\phi$ is weakly differentiable on ${\mathbb R}^d
\backslash {\{} 0 {\}}$ and the weak gradient $\nabla \phi$ (which
is a locally $dx$-integrable function on ${\mathbb R}^d \backslash
{\{} 0 {\}}$), considered as a $dx$-a.e.~defined function on
${\mathbb R}^d$, satisfies
\begin{eqnarray*}
\nabla \phi \in L^1({\mathbb R}^d, \exp(-\phi)  \,dx) \, \cap \,
L^2({\mathbb R}^d, \exp(-\phi)  \,dx).
\end{eqnarray*}
\end{description}
Note that, for many typical potentials in Statistical Physics,
we have $\phi \in
C^{\infty}({\mathbb R}^d\backslash {\{} 0 {\}})$.
For such ``outside the origin regular'' potentials,
condition (D) nevertheless does not exclude a singularity at the point
$0 \in {\mathbb R}^d$.

\begin{description}
\item[(LS)] ({\it local summability})
For all $\Lambda \in {\mathcal O}_c({\mathbb R}^d)$ and all $\gamma \in S_\infty$
\begin{eqnarray*}
\lim_{n\to\infty}\sum_{y\in\gamma
_{\Lambda_n\setminus\Lambda}}\nabla\phi(\cdot- y)
\end{eqnarray*}
exists in $L^1_{{\rm loc}}(\Lambda,\,dx)$.
\end{description}

Assuming $(\phi, \beta, z)$ satisfies (SS), (UI), (LR), (D) and (LS),
and that $\mu$ is the corresponding Gibbs measure constructed with
empty boundary condition, one can prove an integration by parts
formula for the gradient $\nabla^{\Gamma}$, see \cite{AKR98b}, Theorem 4.3.
Utilizing this formula we obtain
for $F, G \in {\mathcal F}C_b^{\infty}({\mathcal D}, \Gamma)$:
\begin{eqnarray}\label{eq99}
{\mathcal E}_{\mu}^{\Gamma}(F,G) = \int_{\Gamma} H^{\Gamma}_{\mu} F G \,d \mu,
\end{eqnarray}
where
\begin{align}\label{eq37}
&H^{\Gamma}_{\mu}F(\gamma) = - \sum_{i,j=1}^{N} \partial_i
\partial_j g_F(\langle f_1, \gamma \rangle, \ldots, \langle f_N,
\gamma \rangle) \langle (\nabla f_i, \nabla f_j)_{{\mathbb R}^d},
\gamma \rangle \\
&- \sum_{j=1}^{N} \partial_j g_F(\langle f_1,
\gamma \rangle, \ldots, \langle f_N, \gamma \rangle)
\Big{(} \langle \Delta f_j, \gamma \rangle  - \beta \sum_{\{x ,y
\} \subset \gamma} \Big(\nabla \phi(x-y), \nabla f_j (x) - \nabla
f_j (y)\Big)_{{\mathbb R}^d} \Big{)}, \nonumber
\end{align}
for $\mu \mbox{-a.e.} \, \gamma \in \Gamma$ and
$F\in {\mathcal F}C_{b}^\infty ({\mathcal D}, \Gamma )$ as in (\ref{eq54}).
Moreover, $H_\mu^\Gamma F \in L^2(\mu)$ for each $F\in {\mathcal
F}C_{b}^\infty ({\mathcal D}, \Gamma )$, see \cite{AKR98b}, Lemma 4.1.
Utilizing (\ref{eq99}), in \cite{AKR98b}, Proposition 5.1, the
following statement has been proven.

\begin{proposition}\label{pr66}
Assume that $(\phi, \beta, z)$ fulfill conditions {\rm (SS)},
{\rm (UI)}, {\rm (LR)}, {\rm (D)}, {\rm (LS)},
and let $\mu$ be the corresponding Gibbs measure constructed with
empty boundary condition. Then the
bilinear form $({\mathcal E}_{\mu}^{\Gamma}, {\mathcal
F}C_b^{\infty}\!({\mathcal D}, \Gamma)\!)$ is closable on $L^2(\mu)$ and
its closure $({\mathcal E}_{\mu}^{\Gamma}, D({\mathcal
E}_{\mu}^{\Gamma})\!)$ is a symmetric Dirichlet form which is
conservative. Its generator is the Friedrichs extension of
$H^{\Gamma}_{\mu}$, which will be denoted by the same symbol.
\end{proposition}

Of course, $H^{\Gamma}_{\mu}$ generates a strongly continuous
contraction semi-group
\begin{eqnarray*}
T^{\mu}_t := \exp(-t H^{\Gamma}_{\mu}), \qquad t \ge 0.
\end{eqnarray*}
The existence of the diffusion process corresponding to
$({\mathcal E}_{\mu}^{\Gamma},
D({\mathcal E}_{\mu}^{\Gamma}))$ was shown in
\cite{AKR98b}, Theorem 5.2, and \cite{MaRo00}, Theorem 4.13.
For all $d \ge 1$ it lives on the
bigger state space $\ddot{\Gamma}$ consisting of all
integer-valued Radon measures on ${\mathbb R}^d$, see
e.g.~\cite{Ka75}. For $d \ge 2$ in \cite{RS98}, Corollary 1, the
authors have proven that the set $\ddot{\Gamma} / \Gamma$ is
${\mathcal E}_{\mu}^{\Gamma}$-exceptional. Thus, the associated
diffusion process can be restricted to a process on $\Gamma$. For
simplicity of notations, we exclude the case $d = 1$ in what
follows. However, all our further considerations do also work in
that case.

\begin{theorem}\label{th57}
Let $(\phi, \beta, z)$ fulfill the same conditions as in
Proposition \ref{pr66} and let $\mu$ be the corresponding Gibbs measure
constructed with empty boundary condition.
Then:

\noindent
{\rm (i)} There exists a conservative diffusion process
(i.e., a conservative strong Markov process with continuous sample paths)
\begin{eqnarray*}
{\bf M} = ({\bf \Omega}, \hat{\bf F}, (\hat{\bf F}_t)_{t \ge 0}, ({\bf \Theta}_t)_{t \ge 0},
({\bf X}(t))_{t \ge 0},
({\bf P}_\gamma)_{\gamma \in {\Gamma}})
\end{eqnarray*}
on ${\Gamma}$ which is properly associated with $({\mathcal E} _{\mu}^{\Gamma},
D({\mathcal E}_{\mu}^{\Gamma}))$, i.e., for all ($\mu$-versions) of $F \in L^2({\Gamma}, \mu)$
and all $t > 0$ the function
\begin{eqnarray*}
\gamma \mapsto p(t, F)(\gamma) := \int_{\bf \Omega} F({\bf X}(t)) \,d{\bf P}_{\gamma},
\qquad \gamma \in {\Gamma},
\end{eqnarray*}
is an ${\mathcal E}_{\mu}^{\Gamma}$-quasi-continuous version of $T^{\mu}_tF$.
The process ${\bf M}$ is up to $\mu$-equivalence unique, has $\mu$ as an invariant measure
and is called microscopic stochastic dynamics.

\noindent
{\rm (ii)} The diffusion process ${\bf M}$ is up to $\mu$-equivalence the unique
diffusion process having $\mu$ as invariant measure and solving
the martingale problem for $(-H^{\Gamma}_{\mu},
D(H^{\Gamma}_{\mu}))$, i.e., for all $G \in D(H^{\Gamma}_{\mu})$
\begin{eqnarray*}
G({\bf X}(t)) - G({\bf X}(0)) + \int_0^t H^{\Gamma}_{\mu} G({\bf
X}(s))\,ds, \qquad t \ge 0,
\end{eqnarray*}
is an $\hat{\bf F}_t$-martingale under ${\bf P}_\gamma$ (hence
starting in $\gamma$) for ${\mathcal E}_{\mu}^{\Gamma}$-q.a.~$\gamma \in
\Gamma$.
\end{theorem}

In the above theorem ${\bf M}$ is canonical, i.e., ${\bf \Omega} =
C([0, \infty) \to {\Gamma}), \, {\bf X}(t)(\xi)$ $= \xi(t), \xi \in {\bf
\Omega}$. The filtration $(\hat{\bf F}_t)_{t \ge 0}$ is the natural
``minimum completed admissible filtration'', cf.~\cite{FOT94}, Chap.~A.2, or
\cite{MaRo92}, Chap.~IV, obtained from
$\sigma{\{}\langle f, {\bf X}(s)\rangle \, | \, 0 \le s \le t,
\, f \in {\mathcal D} {\}}, \, t \ge 0$.
$\hat{\bf F}   :=   \hat{\bf F}_\infty := \bigvee_{t \in [0, \infty)}
\hat{\bf F}_t$
is the smallest $\sigma$-algebra containing all $\hat{\bf F}_t$
and $({\bf \Theta}_t)_{t \ge 0}$ are the corresponding natural time shifts.
For a
detailed discussions of these objects and the notion of quasi-continuity we refer to \cite{MaRo92}.
The second part of the above theorem was proved in
\cite{AKR98b}, Theorem 5.3.

\begin{remark}\label{rm59}
{\rm Let us consider the diffusion process $({\bf X}(t))_{t \ge 0}$
provided by Theorem \ref{th57}. In (\ref{eq37}) we have an
explicit formula for the action of the associated generator
$-H^{\Gamma}_{\mu}$ on smooth cylinder functions. Utilizing an
extension of It{\^o}'s formula to this infinite dimensional
situation on a heuristic level we find the associated infinite system of
stochastic differential equations}:
\begin{align}\label{eq40}
dx(t) & = -\beta \sum_{\stackunder{y(t) \neq x(t)} {y(t) \in {\bf X}(t)}}
\nabla \phi(x(t) - y(t)) \,dt + \sqrt{2} \,dB^x(t),  & x(t) \in {\bf
X}(t), \nonumber \\
{\bf X}(0) & = \gamma,  &\gamma \in \Gamma,
\end{align}
{\rm where $(B^x)_{x = x(0) \in {\bf X}(0)}$ is a sequence of
independent Brownian motions. Theorem \ref{th57}{\rm (ii)} implies that
the process $(({\bf X}(t))_{t \ge 0}, {\bf P}_\gamma)$ solves
the infinite system (\ref{eq40}) in the sense of the associated
martingale problem for ${\mathcal E}_{\mu}^{\Gamma}$-q.a.~$\gamma \in
\Gamma$ as a starting point. \rm}
\end{remark}

\section{Scaling of stochastic dynamics and asso\-ciated Dirichlet form}

We perform the scaling of the process $({\bf X}(t))_{t \ge 0}$ in two steps.

\medskip
\noindent
{\bf First scaling:}
We scale the position of the particles inside the configuration space as follows:
\begin{eqnarray*}
\Gamma \ni \gamma \mapsto S_{in, \epsilon}(\gamma)
:= \{ \epsilon \, x | \, x \in \gamma \} \in \Gamma, \qquad \epsilon > 0,
\end{eqnarray*}
i.e., for $f \in {\mathcal D}$, the scaling is given through
$\langle f, S_{in, \epsilon}(\gamma) \rangle
= \sum_{x \in \gamma} f(\epsilon \, x)$.
Obviously, $S_{in, \epsilon}$ is a homeomorphism on $\Gamma$.
From now on we assume that $\mu$ corresponds to $(\phi, \beta, 1), \, \beta \ge 0$
and the construction with empty boundary condition. Let us define the image measure
$\tilde{\mu}_\epsilon :=  S_{in, \epsilon}^{\ast} \mu$.
This measure is also defined on $(\Gamma, {\mathcal B}(\Gamma))$ and
it is easy to check that it is the Gibbs measure corresponding to $(\phi_\epsilon,
\beta, \epsilon^{-d})$ and the construction with empty boundary condition,
where $\phi_{\epsilon} := \phi(\epsilon^{-1} \cdot)$.
Furthermore, since
$C(\beta \phi_\epsilon, \epsilon^{-d}) = C(\beta \phi, 1)$,
recall (UI), the measure $\tilde{\mu}_{\epsilon}$ is
in the LA-HT regime if and only if this is true for $\mu$.

\medskip
\noindent
{\bf Second scaling:}
This scaling leads us out of the configuration space and is given by
\begin{eqnarray*}
\Gamma \ni \gamma \mapsto S_{out, \epsilon}(\gamma) :=
\epsilon^{d/2}\Big(\gamma - \rho^{(1)}_{\phi_\epsilon}(\beta,
\epsilon^{-d}) \, \epsilon^{-d} \,dx \Big) \in \Gamma_\epsilon
\end{eqnarray*}
where $\Gamma_\epsilon :=  S_{out, \epsilon}(\Gamma) \subset {\mathcal
D}^{\prime}, \, \epsilon > 0$, ${\mathcal D}^\prime$ is the
topological dual of ${\mathcal D}$ (where both ${\mathcal D}$ and ${\mathcal
D}^\prime$ are equipped with their respective usual locally convex topology).
We consider $\Gamma_{\epsilon}$ as a topological subspace of
${\mathcal D}^\prime$, thus $\Gamma_{\epsilon}$ is equipped with the
corresponding Borel $\sigma$-algebra. Obviously, $S_{out, \epsilon}: \Gamma
\to \Gamma_{\epsilon}$ is continuous, hence Borel-measurable. Since it is
also one-to-one and since both $\Gamma$ and ${\mathcal D}^\prime$ are standard
measurable spaces, it follows by \cite{Par67}, Chap.~V, Theorem 2.4,
that $\Gamma_\epsilon$ is a
Borel subset of ${\mathcal D}^\prime$ and that $S^{-1}_{out, \epsilon}:
\Gamma_{\epsilon} \to \Gamma$ is also Borel-measurable.
The function
$\rho^{(1)}_{\phi_\epsilon}(\beta, \epsilon^{-d})$ is the first
correlation function corresponding to the Gibbs measure
$\tilde{\mu}_\epsilon$, i.e.,
\begin{eqnarray*}
\int_{{\mathbb R}^d} f(x) \, \rho^{(1)}_{\phi_\epsilon}(\beta,
\epsilon^{-d}) \, \epsilon^{-d} \,dx
= \int_{\Gamma} \langle f, \gamma \rangle \,d\tilde{\mu}_\epsilon(\gamma),
\qquad \forall f \in C_0({\mathbb R}^d).
\end{eqnarray*}
Applied to a test function $f \in {\mathcal
D}$, the second scaling gives
\begin{eqnarray}\label{eq1972}
\langle f, S_{out, \epsilon}(\gamma) \rangle = \epsilon^{d/2}
\Big{(} \sum_{x \in \gamma}  f(x) -
\rho^{(1)}_{\phi_\epsilon}(\beta, \epsilon^{-d}) \, \epsilon^{-d}
\int f(x) \,dx \Big{)},
\end{eqnarray}
where $\langle \cdot, \cdot \rangle$ denotes the dual paring
between ${\mathcal D}$ and ${\mathcal D}^{\prime}$. Here we assume the LA-HT
regime. So, as mentioned before $\rho^{(1)}_{\phi}(\beta,1)$ is a
constant, and thus by definition of
$\tilde{\mu}_\epsilon$
also $\rho^{(1)}_{\phi_\epsilon}(\beta,
\epsilon^{-d})$ is a constant, see Subsection \ref{ss23}. Obviously,
the random variable (\ref{eq1972}) is centered
w.r.t.~the measure $\tilde{\mu}_\epsilon$.

\medskip
\noindent
{\bf Scaled process:}
The scaled process of our interest is
\begin{eqnarray*}
{\bf X}_{\epsilon}(t) := S_{out, \epsilon} (S_{in, \epsilon} ({\bf X}
(\epsilon^{-2}t))), \qquad t \ge 0, \quad \epsilon > 0.
\end{eqnarray*}

\medskip
\noindent
{\bf Associated Dirichlet form:}

Next for each $\epsilon > 0$ we construct a Dirichlet form ${\mathcal E}_{\epsilon}$ such that
$({\bf X}_{\epsilon}(t))_{t \ge 0}$ is the unique process which is properly associated
to ${\mathcal E}_{\epsilon}$.

Let $\mu_\epsilon :=  S_{out, \epsilon}^{\ast} S_{in,
\epsilon}^{\ast} \mu = S_{out, \epsilon}^{\ast}
\tilde{\mu}_\epsilon$. Then we define a unitary mapping ${\mathcal
S}_{out, \epsilon}:  L^2(\Gamma_{\epsilon}, \mu_\epsilon) \to
L^2(\Gamma, \tilde{\mu}_\epsilon)$ by defining ${\mathcal S}_{out,
\epsilon}F$ to be the $\tilde{\mu}_\epsilon$-class represented by
$\tilde{F} \circ S_{out, \epsilon}$ for any
$\mu_{\epsilon}$-version $\tilde{F}$ of $F \in
L^2(\Gamma_{\epsilon}, \mu_\epsilon)$. Using this mapping we
define a bilinear form $({\mathcal E}_{\epsilon}, D({\mathcal
E}_{\epsilon}))$ as the image bilinear form of $({\mathcal
E}_{\tilde{\mu}_\epsilon}^{\Gamma}, D({\mathcal
E}_{\tilde{\mu}_\epsilon}^{\Gamma}))$ under the mapping ${\mathcal
S}_{out, \epsilon}$:
\begin{eqnarray}\label{eq124}
{\mathcal E}_{\epsilon}(F,G) := {\mathcal E}_{\tilde{\mu}_\epsilon}^{\Gamma}
({\mathcal S}_{out, \epsilon} F, {\mathcal S}_{out, \epsilon} G), \qquad F, G \in D({\mathcal
E}_{\epsilon}),
\end{eqnarray}
where
$D({\mathcal E}_{\epsilon}) := {\mathcal S}^{-1}_{out, \epsilon} D({\mathcal
E}_{\tilde{\mu}_\epsilon}^{\Gamma})$.
Let ${\mathcal F}C_b^{\infty}({\mathcal D}, \Gamma_{\epsilon})$ be defined
analogously to the space ${\mathcal F}C_b^{\infty}({\mathcal D}, \Gamma)$.
Then obviously, ${\mathcal F}C_b^{\infty}({\mathcal D}, \Gamma_{\epsilon})
\subset D({\mathcal E}_{\epsilon})$, hence $({\mathcal
E}_{\epsilon}, D({\mathcal
E}_{\epsilon}))$ is densely defined. It follows by
\cite{MaRo92}, Chapter VI, Exercise 1.1, that
$({\mathcal E}_{\epsilon}, D({\mathcal E}_{\epsilon}))$ is a Dirichlet form.
It is called the image Dirichlet form of $({\mathcal
E}_{\tilde{\mu}_\epsilon}^{\Gamma}, D({\mathcal
E}_{\tilde{\mu}_\epsilon}^{\Gamma}))$ under the mapping $S_{out,
\epsilon}$. Its generator $(-H_\epsilon, D(H_\epsilon))$ is given
by
\begin{eqnarray}\label{eq25}
H_\epsilon = {\mathcal S}^{-1}_{out, \epsilon}
H^{\Gamma}_{\tilde{\mu}_\epsilon} {\mathcal S}_{out, \epsilon}, \qquad
D(H_{\epsilon}) =  {\mathcal S}^{-1}_{out, \epsilon}
D(H_{\tilde{\mu}_\epsilon}^{\Gamma}).
\end{eqnarray}
Then for $F \in {\mathcal F}C_b^{\infty}({\mathcal D}, \Gamma_{\epsilon})
\subset D(H_\epsilon)$ we have
\begin{gather}\label{eq1147}
H_{\epsilon} F(\omega) =  - \sum_{i,j=1}^{N} \frac{\partial^2
g_F}{\partial s_i \partial s_j} (\langle f_1, \omega \rangle,
\ldots, \langle f_N, \omega \rangle) \\ \times \Big\langle
(\nabla f_i, \nabla f_j)_{{\mathbb R}^d}, \epsilon^{d/2}\omega +
\rho^{(1)}_{\phi_\epsilon}(\beta, \epsilon^{-d}) \,dx
\Big\rangle - \sum_{j=1}^{N} \partial_j g_F
(\langle f_1, \omega \rangle, \ldots, \langle f_N, \omega \rangle)
\nonumber \\
\times \Big{(} \langle \Delta f_j, \omega
\rangle  - \epsilon^{d/2} \beta \sum_{\{x ,y \} \subset
S^{-1}_{out, \epsilon}\omega} \Big(\nabla \phi_{\epsilon}(x-y),
\nabla f_j (x) - \nabla f_j (y)\Big)_{{\mathbb R}^d} \Big{)},
\nonumber
\end{gather}
where $F$ is of the form (\ref{eq54}) and the variable $\omega$ is
running through $\Gamma_\epsilon$. Note that the last term is well-defined
for $\mu_{\epsilon}$-a.e.~$\omega \in \Gamma_\epsilon$.

\begin{theorem}\label{th58}
Let $(\phi, \beta, 1)$ fulfill conditions {\rm (SS)}, {\rm (UI)}, {\rm (LR)},
{\rm (D)}, {\rm (LS)} and $\mu$ be the corresponding Gibbs measure
constructed with empty boundary condition. Then for
all ($\mu_{\epsilon}$-versions) of $F \in L^2({\Gamma}_\epsilon,
\mu_\epsilon)$ and all $t > 0$ the function
\begin{eqnarray*}
\omega \mapsto p_\epsilon(t, F)(\omega) :=
\int_{{\bf \Omega}} F({\bf X}_{\epsilon}(t))
\,d{\bf P}_{S^{-1}_{in, \epsilon} S^{-1}_{out, \epsilon} \omega},
\qquad \omega \in {\Gamma}_\epsilon,
\end{eqnarray*}
is a $\mu_{\epsilon}$-version of $T_{\epsilon, t}F := \exp(-t H_\epsilon)F$.
For ${\bf Q}_{\omega} := {\bf P}_{S^{-1}_{in, \epsilon}
S^{-1}_{out, \epsilon} \omega}, \, \omega \in \Gamma_{\epsilon}$, the process
${\bf M}^{\epsilon} = ({\bf \Omega}, \hat{\bf F},
(\hat{\bf F}_{t/\epsilon^2})_{t \ge 0},
({\bf \Theta}_{t/\epsilon^2})_{t \ge 0}, ({\bf X}_{\epsilon}(t))_{t \ge 0},
({\bf Q}_\omega)_{\omega \in {\Gamma}_\epsilon})$ is a diffusion process and
thus up to $\mu_\epsilon$-equivalence the unique process in this class which is
properly associated with $({\mathcal E}_{\epsilon}, D({\mathcal E}_{\epsilon}))$
and has $\mu_\epsilon$ as an invariant measure.
\end{theorem}
{\bf Proof:}
For $F \in L^2(\Gamma_{\epsilon}, \mu_\epsilon)$ we have
$F({\bf X}_{\epsilon}(t)) = ({\mathcal S}_{in, \epsilon}{\mathcal S}_{out, \epsilon}
F)({\bf X}(\epsilon^{-2}t)), \, t \ge 0$,
where ${\mathcal S}_{in, \epsilon}F := F \circ S_{in, \epsilon}$.
By Theorem
\ref{th57} we have
\begin{multline}\label{eq81}
({\mathcal S}^{-1}_{out, \epsilon} {\mathcal S}^{-1}_{in, \epsilon} \exp(-t \epsilon^{-2}
H^{\Gamma}_{\mu}) {\mathcal S}_{in, \epsilon}{\mathcal S}_{out, \epsilon}F)(\omega) \\
= \int_{{\bf \Omega}} {\mathcal S}_{out, \epsilon}{\mathcal S}_{in, \epsilon}F({\bf
X}(\epsilon^{-2}t)) \,d{\bf P}_{S^{-1}_{in, \epsilon} S^{-1}_{out,
\epsilon} \omega} = \int_{{\bf \Omega}} F({\bf
X}_{\epsilon}(t)) \,d{\bf P}_{S^{-1}_{in, \epsilon} S^{-1}_{out,
\epsilon} \omega},
\end{multline}
for $\mu_{\epsilon}$ almost all $\omega \in \Gamma_{\epsilon}$.
We note that $({\mathcal E}^{\Gamma}_{\tilde{\mu_\epsilon}},
D({\mathcal E}^{\Gamma}_{\tilde{\mu_\epsilon}}))$ is obviously the image
Dirichlet form under the map ${\mathcal S}_{in, \epsilon}$ of
$({\mathcal E}^{\Gamma}_{\mu}, D({\mathcal E}^{\Gamma}_{\mu}))$ times
$\epsilon^{-2}$. Hence we have for the corresponding generator
$(H^{\Gamma}_{\tilde{\mu_\epsilon}}, D(H^{\Gamma}_{\tilde{\mu_\epsilon}}))$
\begin{eqnarray}\label{eq26}
H^{\Gamma}_{\tilde{\mu}_\epsilon} = {\mathcal S}^{-1}_{in, \epsilon}
\epsilon^{-2} H^{\Gamma}_{\mu} {\mathcal S}_{in, \epsilon}, \quad
D(H^{\Gamma}_{\tilde{\mu_\epsilon}})
= {\mathcal S}^{-1}_{in, \epsilon} (D(H_{\mu}^{\Gamma})).
\end{eqnarray}
Using the Hille-Yosida theorem (via resolvent) and (\ref{eq25}),
(\ref{eq26}), we can conclude that
\begin{eqnarray}\label{eq82}
{\mathcal S}^{-1}_{out, \epsilon} {\mathcal S}^{-1}_{in, \epsilon} \exp(-t \epsilon^{-2}
H^{\Gamma}_{\mu}) {\mathcal S}_{in, \epsilon}{\mathcal S}_{out, \epsilon}
= \exp(-t H_\epsilon)
\end{eqnarray}
on $L^2(\Gamma_\epsilon, \mu_\epsilon)$. Thus, by (\ref{eq81}) and
(\ref{eq82}) the first statement of the theorem is proved.
The fact that ${\bf M}^{\epsilon}$ is a diffusion is straightforward to
check. In particular, it then follows by \cite{MaRo92}, Chap.~IV,
Theorem 3.5, that ${\bf M}^{\epsilon}$ is properly associated with
$({\mathcal E}_\epsilon, D({\mathcal E}_\epsilon))$.
\hfill $\blacksquare$

\section{Convergence of Dirichlet forms}

Our aim is to show convergence of $({\bf X}_{\epsilon}(t))_{t \ge
0}$ to a generalized Ornstein--Uhlenbeck process $({\bf X}(t))_{t
\ge 0}$ as $\epsilon \to 0$. In this section
we prove this in terms of the corresponding Dirichlet forms.

It will turn out that the limit Dirichlet form is defined in
$L^2({\mathcal D}^\prime, \nu_{\mu})$, where $\nu_{\mu}$ is the
Gaussian white noise measure on ${\mathcal D}^\prime$ with covariance
operator $\chi_{\phi}(\beta) \operatorname{Id}$ and
\begin{eqnarray*}
\chi_{\phi}(\beta) := \rho_{\phi}^{(1)}(\beta, 1) +
\int_{{\mathbb R}^{d}} u_{\phi}^{(2)}(\beta, 1, x, 0) \,dx
\end{eqnarray*}
is the {\it compressibility} of the Gibbs state $\mu$,
see (\ref{eq22a}) below for the definition of the Ursell function
$u_{\phi}^{(2)}$ and Proposition \ref{pr1969} for the existence
of the integral.
The measure $\nu_{\mu}$ exists due to
the Bochner--Minlos theorem via its characteristic function given by
\begin{eqnarray*}
\int_{{\mathcal D}^\prime} \exp(i\langle f,\omega \rangle) \,d\nu_{\mu}(\omega)
= \exp\Big{(} -\frac{\chi_{\phi}(\beta) }{2}
\int_{{\mathbb R}^{d}} (f(x))^2 \,dx \Big{)}, \quad f \in {\mathcal D}.
\end{eqnarray*}

For $n \in {\mathbb Z}$ we define a weighted Sobolev spaces ${\mathcal H}_n$ as the
closure of ${\mathcal D}$ w.r.t.~the Hilbert norm
\begin{eqnarray*}
\parallel f \parallel^2_{n} \, = \, \langle f,f \rangle_{n}
\, := \int_{{\mathbb R}^d} A^n f(x) f(x) \,dx,
\qquad f \in {\mathcal D},
\end{eqnarray*}
where
$Af(x) = - \Delta f(x) + |x|^2 f(x), \, x \in {\mathbb R}^d$, i.e., $A$
is the Hamilton operator of the harmonic oscillator
with ground state eigenvalue $d$. We identify ${\mathcal H}_0 = L^2({\mathbb R}^d,
dx)$ with its dual and obtain
\begin{eqnarray*}
{\mathcal D} \subset S({\mathbb R}^d) \subset {\mathcal H}_n \subset
L^2({\mathbb R}^d, dx) \subset {\mathcal H}_{-n} \subset  S'({\mathbb
R}^d) \subset {\mathcal D}^{\prime}, \quad n \in {\mathbb N}.
\end{eqnarray*}
Here as usual $S'({\mathbb R}^d)$ denotes the space of tempered
distributions which is the topological dual of $S({\mathbb R}^d)$,
the Schwartz space of smooth functions on
${\mathbb R}^d$ decaying faster than any polynomial.
Of course, ${\mathcal H}_{-n}$ is the topological dual of ${\mathcal H}_{n}$
w.r.t.~${\mathcal H}_{0}$. The dual paring between these spaces we denote by
$\langle \cdot, \cdot \rangle$. Since the embeddings
${\mathcal H}_n \subset {\mathcal H}_{n-d}$ are Hilbert--Schmidt for all
$n \in {\mathbb Z}$, it follows by the Bochner--Minlos theorem that
$\nu_{\mu}({\mathcal H}_{-d}) = 1$.

The first part of the following theorem is an easy generalization of
Proposition 3.9 in \cite{Br80}. The second and third
part have been proved in \cite{Br80}, Proposition 5.4 and Theorem 6.5,
respectively.
\begin{theorem}\label{th9}
Let us assume that $(\phi, \beta, 1)$ fulfill {\rm (S)}, {\rm
(UI)} and let $\mu$ be the corresponding Gibbs measure constructed
with empty boundary condition. Then:

\noindent
{\rm (i)} There exists $C^{(1)} \in (0, \infty)$ such that
\begin{eqnarray*}
\int_{{\mathcal D}^\prime} \parallel \omega \parallel^2_{-(d+1)}
\,d\mu_\epsilon(\omega) \le C^{(1)}
\end{eqnarray*}
uniformly in $\epsilon \in (0,1]$ and, in particular,
$\mu_\epsilon({\mathcal H}_{-(d+1)}) = 1$.

\noindent
{\rm (ii)}
For each $f \in {\mathcal D}$ we have $\lim_{\epsilon \to 0}
{\mathbb E}_{\mu_\epsilon}[\langle f, \cdot \rangle^2]
= {\mathbb E}_{\nu_{\mu}}[\langle f, \cdot \rangle^2]$.

\noindent
{\rm (iii)} The family of measures $({\mu}_{\epsilon})_{\epsilon > 0}$
converges weakly on ${\mathcal H}_{-(d+1)}$ to the Gaussian measure
$\nu_{\mu}$ as $\epsilon \to 0$.
\end{theorem}

We shall also use the following lemma, which is easy to derive by using the
properties of correlation functions, see Section \ref{ss23}, and recalling
that $\tilde{\mu}_\epsilon = S_{in, \epsilon}^{\ast} \mu$ is the Gibbs
measure corresponding to $(\phi_{\epsilon}, \beta, \epsilon^{-d})$ and
the construction with empty boundary condition.
\begin{lemma}\label{le7}
Let the conditions of Theorem \ref{th9} hold. Then we have:
\begin{align*}
\rho^{(1)}_{\phi_\epsilon}(\beta, \epsilon^{-d}) & =
\rho^{(1)}_{\phi}(\beta, 1), \\
\rho^{(2)}_{\phi_\epsilon}(\beta , \epsilon^{-d}, x, y) & =
\rho^{(2)}_{\phi}\Big{(}\beta, 1, \frac{x-y}{\epsilon}, 0 \Big{)}.
\end{align*}
\end{lemma}

We define the Dirichlet form $({\mathcal E}_{\nu_{\mu}},
D({\mathcal E}_{\nu_{\mu}}))$ as the
closure of the bilinear form
\begin{multline*}
{\mathcal E}_{\nu_{\mu}}(F,G) = - \rho_{\phi}^{(1)}(\beta,1) \int_{{\mathcal D}^\prime} \int_{{\mathbb R}^d}
\partial_x F(\omega) \Delta \partial_x G(\omega)
\,dx \,d\nu_{\mu}(\omega) \\
= \rho_{\phi}^{(1)}(\beta,1) \int_{{\mathcal D}^\prime} \int_{{\mathbb R}^d}
(\nabla \partial_x F(\omega), \nabla \partial_x G(\omega))_{{\mathbb R}^d}
\,dx \,d\nu_{\mu}(\omega),
\end{multline*}
where $F, G \in {\mathcal F}C_b^{\infty}({\mathcal D}, {\mathcal D}^{\prime})$ and
the space ${\mathcal F}C_b^{\infty}({\mathcal D}, {\mathcal D}^{\prime})$ is defined
analogously to ${\mathcal F}C_b^{\infty}({\mathcal D}, \Gamma)$. Here
$\partial_x F$ denotes the derivative of $F \! = g_F(\langle f_1, \cdot \rangle,
\ldots, \langle f_N, \cdot \rangle)$
$\in {\mathcal F}C_b^{\infty}({\mathcal D}, {\mathcal D}^{\prime})$ in direction
$\varepsilon_x, \, x \in {\mathbb R}^d$, i.e.,
\begin{align*}
\partial_x F(\omega) = \frac{d}{dt}F(\omega +
t\varepsilon_x)\Big{|}_{t=0}
= \sum_{j=1}^{N} \partial_j g_F
(\langle f_1, \omega \rangle, \ldots,
\langle f_N, \omega \rangle) f_j(x), \quad \omega \in {\mathcal D}^{\prime},
\end{align*}
where $N \in {\mathbb N}$ and $f_1, \ldots, f_N \in {\mathcal D}$.

Integrating by parts in the Gaussian space, see e.g.~\cite{BeKo88},
Theorem 6.1.2 and 6.1.3, we obtain
\begin{eqnarray*}
{\mathcal E}_{\nu_{\mu}}(F,G) = \int_{{\mathcal D}^{\prime}}
H F(\omega) G(\omega) \,d\nu_{\mu}(\omega),
\qquad F, G \in {\mathcal F}C_b^{\infty}({\mathcal D}, {\mathcal D}^{\prime}),
\end{eqnarray*}
where
\begin{align}\label{eq5}
H F & = - \rho_{\phi}^{(1)}(\beta,1) \sum_{i,j=1}^{N}
\partial_i \partial_j g_F (\langle f_1,  \cdot \rangle,
\ldots, \langle f_N,  \cdot \rangle) \int_{{\mathbb R}^d} (\nabla
f_i(x), \nabla f_j(x))_{{\mathbb R}^d} \,dx \nonumber \\ & -
\frac{\rho_{\phi}^{(1)}(\beta,1)}{\chi_{\phi}(\beta) }
\sum_{j=1}^{N} \partial_j g_F (\langle f_1,  \cdot \rangle,
\ldots, \langle f_N,  \cdot \rangle) \langle \Delta f_j,  \cdot
\rangle.
\end{align}
It is well-known, see e.g.~\cite{BeKo88}, Theorem 6.1.4, that the
operator $H$ is essentially self-adjoint on ${\mathcal
F}C_b^{\infty}({\mathcal D}, {\mathcal D}^{\prime})$. We preserve the same
notation for its closure. The operator $H$ generates an infinite
dimensional Ornstein--Uhlenbeck semi-group
\begin{eqnarray*}
T_t := \exp(-tH), \qquad t \ge 0,
\end{eqnarray*}
in $L^2(\nu_\mu)$. This semi-group is associated to a generalized
Ornstein--Uhlenbeck process $({\bf X}(t))_{t \ge 0}$ on ${\mathcal
D}^{\prime}$, see \cite{BeKo88}, Chapter 6, Section 1.5.

\begin{theorem}\label{th11}
Suppose that $(\phi, \beta, 1)$ satisfy the conditions {\rm (S)},
{\rm (UI)} and let $\mu$ be the corresponding Gibbs measure
construction with empty boundary condition.
Then for all $F, G \in {\mathcal F}C_b^{\infty}({\mathcal D},
{\mathcal D}^{\prime})$ we have
\begin{eqnarray}\label{eq357}
\lim_{\epsilon \to 0} {\mathcal E}_{\epsilon}(F,G) = {\mathcal E}_{\nu_{\mu}}(F,G).
\end{eqnarray}
\end{theorem}

\begin{remark}\label{rm1111}
{\rm (i) The $({\bf X}(t))_{t \ge 0}$ is the unique process
associated to the closure of the pre-Dirichlet form $({\mathcal E}, {\mathcal
F}C_b^{\infty}({\mathcal D}, {\mathcal D}^{\prime}))$ on
$L^2({\mathcal D}^{\prime}, \nu_\mu)$. In this sense the
convergence of bilinear forms proven in Theorem \ref{th11}
uniquely determines the limiting process $({\bf X}(t))_{t \ge 0}$.

\noindent
(ii) The generator $H$ corresponds to the following
stochastic differential equation:
\begin{eqnarray*}
d{\bf X}(t,x) =
\frac{\rho_{\phi}^{(1)}(\beta,1)}{\chi_{\phi}(\beta) } \Delta {\bf
X}(t,x) \,dt + \sqrt{2 \, \rho_{\phi}^{(1)}(\beta,1)} \,d{\bf W}(t,x),
\end{eqnarray*}
where $({\bf W}(t))_{t \ge 0}$ is a Brownian motion in ${\mathcal
D}^{\prime}$ with covariance operator $- \Delta$, and the
coefficient $\rho_{\phi}^{(1)}(\beta,1) / \chi_{\phi}(\beta) $ is
called {\it bulk diffusion coefficient}.

\noindent
(iii) The generality of the class of admissible potentials is very important
from the physical point of view. Before one could only treat smooth,
compactly supported, positive potentials. However, any physical realistic
potential has a singularity at the origin. Furthermore, it is of physical
interest to study potentials which also have a negative part.

\noindent
(iv) The proof of Theorem \ref{th11} is straightforward. However, it identifies
the bulk diffusion coefficient for very general potentials. This
coefficient is, from the physical point of view, the most interesting quantity.}
\end{remark}
{\bf Proof:} We first note that each function
$F \in {\mathcal F}C_b^{\infty}({\mathcal D}, {\mathcal D}^{\prime})$,
when restricted to $\Gamma_\epsilon$, belongs to ${\mathcal
F}C_b^{\infty}({\mathcal D}, \Gamma_\epsilon)$ $\subset D({\mathcal
E}_{\epsilon})$. Furthermore, since ${\mathcal B}({\mathcal D}^\prime)
\cap \Gamma_\epsilon = {\mathcal B}(\Gamma_\epsilon)$, the measure
$\mu_\epsilon$ can be considered as a measure on $({\mathcal
D}^\prime, {\mathcal B}({\mathcal D}^\prime))$. By the
polarization identity, it is sufficient to prove (\ref{eq357}) for
the case $G = F = g_F(\langle f_1, \omega \rangle,
\ldots, \langle f_N, \omega \rangle)$.
Evaluating (\ref{eq124}) and applying Lemma \ref{le7}
we obtain
\begin{multline}\label{eq50}
 {\mathcal E}_{\epsilon}(F,F) = \epsilon^{d} \sum_{i,j=1}^{N} \int_{\Gamma}
\langle (\nabla f_i, \nabla f_j)_{{\mathbb R}^d}, \gamma \rangle
\partial_i g_F
(\langle f_1,\epsilon^{d/2}(\gamma - \epsilon^{-d} \rho_{\phi_\epsilon}^{(1)}(\beta,\epsilon^{-d}) \,dx) \rangle, \\
\ldots,
\langle f_N, \epsilon^{d/2}(\gamma - \epsilon^{-d} \rho_{\phi_\epsilon}^{(1)}(\beta,\epsilon^{-d}) \,dx) \rangle)
\partial_j g_F
(\langle f_1, \epsilon^{d/2}(\gamma - \epsilon^{-d} \rho_{\phi_\epsilon}^{(1)}(\beta,\epsilon^{-d})  \,dx) \rangle, \\ \ldots,
\langle f_N, \epsilon^{d/2}(\gamma - \epsilon^{-d} \rho_{\phi_\epsilon}^{(1)}(\beta,\epsilon^{-d})  \,dx) \rangle)
\,d\tilde{\mu}_\epsilon(\gamma) \\
= \epsilon^{d/2} \sum_{i,j=1}^{N} \int_{{\mathcal D}^{\prime}}
\langle (\nabla f_i, \nabla f_j)_{{\mathbb R}^d}, \omega \rangle
\partial_i g_F
(\langle f_1,\omega \rangle, \ldots,
\langle f_N, \omega \rangle) \partial_j g_F
(\langle f_1, \omega \rangle, \ldots,
\langle f_N, \omega \rangle)
\,d\mu_\epsilon(\omega) \\
+
\rho_{\phi}^{(1)}(\beta,1) \sum_{i,j=1}^{N} \int_{{\mathbb R}^d}
(\nabla f_i(x), \nabla f_j(x))_{{\mathbb R}^d} \,dx \\
\times \int_{{\mathcal D}^{\prime}} \partial_i g_F
(\langle f_1,\omega \rangle, \ldots,
\langle f_N, \omega \rangle) \partial_j g_F
(\langle f_1, \omega \rangle, \ldots,
\langle f_N, \omega \rangle)
\,d\mu_\epsilon(\omega).
\end{multline}
By Theorem \ref{th9}(iii) we get
\begin{gather*}
\lim_{\epsilon \to 0} \int_{{\mathcal D}^{\prime}} \partial_i g_F
(\langle f_1,\omega \rangle, \ldots,
\langle f_N, \omega \rangle) \partial_j g_F
(\langle f_1, \omega \rangle, \ldots,
\langle f_N, \omega \rangle)
\,d\mu_\epsilon(\omega) \\
=  \int_{{\mathcal D}^{\prime}} \partial_i g_F
(\langle f_1,\omega \rangle, \ldots,
\langle f_N, \omega \rangle) \partial_j g_F
(\langle f_1, \omega \rangle, \ldots,
\langle f_N, \omega \rangle) \nu_{\mu}(\omega),
\end{gather*}
hence, the second term in (\ref{eq50}) converges to ${\mathcal
E}_{\nu_{\mu}}(F,F)$ and it only remains to show that first term
in (\ref{eq50}) converges to zero as $\epsilon \to 0$.
But this is obvious from Theorem \ref{th9}(i), because
$F \in {\mathcal F}C_b^{\infty}({\mathcal D}, {\mathcal D}^{\prime})$.
\hfill $\blacksquare$

\section{Convergence in law}

The convergence in terms of the Dirichlet forms admits no
probabilistic interpretation. Hence, next we study convergence in
law of the scaled processes.

The laws of the scaled equilibrium processes
${\bf P}^{\epsilon} :=  {\bf Q}_{\mu_\epsilon} \circ {\bf X}_{\epsilon}^{-1}
(= {\bf P}_{\mu} \circ {\bf X}_{\epsilon}^{-1})$,
are probability measures on $C([0, \infty), \Gamma_\epsilon)$, where
${\bf Q}_{\mu_\epsilon} := \int_{\Gamma_\epsilon} {\bf Q}_\omega
\,d\mu_\epsilon(\omega)$ and ${\bf P}_{\mu} := \int_{\Gamma} {\bf P}_\gamma
\,d\mu (\gamma)$,
cf.~Theorem \ref{th58}. Since $C([0, \infty), \Gamma_\epsilon)$ is a Borel subset
of $C([0, \infty), {\mathcal D}^{\prime})$ (under the natural embedding)
with compatible measurable structures we can consider
${\bf P}^{\epsilon}$ as a measure on $C([0, \infty), {\mathcal D}^{\prime})$
and by using Theorem \ref{th57}{\rm (ii)}
we find that the process $({\bf X}(t))_{t \ge 0}$ corresponding to
${\bf P}^{\epsilon}$, i.e., the realization of $({\bf
X}_{\epsilon}(t))_{t \ge 0}$ as a coordinate process in $C([0,
\infty), {\mathcal D}^{\prime})$, solves the martingale problem for
$(-H_{\epsilon}, D(H^{\Gamma}_{\mu}))$ w.r.t.~the corresponding
minimum completed admissible filtration
$({\bf F}_t)_{t \ge 0}$ for all $\epsilon > 0$.

\subsection{Tightness}

\begin{theorem}\label{pr777}
Let $(\phi, \beta, 1)$ fulfill conditions {\rm (SS)}, {\rm
(UI)}, {\rm (LR)}, {\rm (D)}, {\rm (LS)} and let $\mu$
be the corresponding Gibbs measure constructed with
empty boundary condition. Then there exists $m
\in {\mathbb N}, \, m \ge d+1,$ such that the family of probability measures
$({\bf P}^\epsilon)_{\epsilon > 0}$ can be restricted to
the space $C([0, \infty), {\mathcal H}_{-m})$. Furthermore, $({\bf
P}^\epsilon)_{\epsilon > 0}$ is tight on $C([0, \infty),
{\mathcal H}_{-m})$.
\end{theorem}
\begin{remark}\label{rm0002}
{\rm Theorem \ref{pr777} has been proved before by
T.~Brox, \cite{Br80}, and H.~Spohn, \cite{Sp86}, for smooth, compactly
supported potentials only. Their proof can be generalized to a more general class of
potentials. However, their technique requires that $\partial_j\phi x^i$
(here $\partial_j\phi$ is the partial derivative of the potential in direction
$j$ and $x^i$ the $i$-th component of the identity)
is locally integrable w.r.t.~the Lebesgue measure. From the physical point of
view this is a very restrictive assumption on the singularity of
the potential at the origin.}
\end{remark}
{\bf Proof:}
Let $f \in {\mathcal D}$. By
Theorem \ref{th9}(i) we know, in particular, that the functions
$\langle f, \cdot \rangle$, $\langle \nabla f, \cdot \rangle \in
L^2(\mu_\epsilon)$. Hence it is easy to show by approximation that
$\langle f, \cdot \rangle \in D({\mathcal E}_{\epsilon})$.
Consider the conservative diffusion process ${\bf M}^{\epsilon}$ on
$\Gamma_\epsilon$ associated with
$({\mathcal E}_{\epsilon}, D({\mathcal E}_{\epsilon}))$ according to
Theorem \ref{th58}. We may regard ${\bf M}^{\epsilon}$ on the
state space ${\mathcal D}^{\prime}$ (common to all
${\bf M}^{\epsilon}, \epsilon > 0$). Considering its distribution
on $C([0, \infty), {\mathcal D}^{\prime})$ we may regard its canonical
realization ${\bf M}^{\epsilon} = ({\bf \Omega}, {\bf F},
({\bf F}_{t})_{t \ge 0}, ({\bf \Theta}_{t})_{t \ge 0},
({\bf X}(t))_{t \ge 0}, ({\bf Q}^{\epsilon}_\omega)_{\omega
\in {\mathcal D}^{\prime}})$. So, in particular,
$\Omega = C([0, \infty), {\mathcal D}^{\prime}), \, X(t)(\omega) =
\omega(t), t \ge 0, \, \theta_t(\omega) = \omega(t + \cdot)$, and
${\bf P}^{\epsilon} = \int_{\Gamma_\epsilon} {\bf Q}^{\epsilon}_\omega
\,d\mu_\epsilon(\omega)$. Fix $T > 0$. Below we canonically
project the process onto $\Omega_T := C([0, T], {\mathcal D}^{\prime})$
without expressing this explicitly. We define the time reversal
$r_{\scriptscriptstyle{T}}(\omega) := \omega(T- \cdot),
\, \omega \in \Omega_T$. Now, by the well-known Lyons-Zheng
decomposition, cf. \cite{LZ88}, \cite{FOT94}, and also \cite{LZ94}
for its infinite dimensional variant, we have for all $0 \le t \le T$:
\begin{eqnarray*}
\langle f, {\bf  X}(t) \rangle - \langle f, {\bf  X}(0) \rangle
= \frac{1}{2} {\bf  M}_t(\epsilon, f)
+ \frac{1}{2}\Big{(}{\bf M}_{T-t}(\epsilon, f)(r_{\scriptscriptstyle{T}})
- {\bf  M}_T(\epsilon, f)(r_{\scriptscriptstyle{T}}) \Big{)}
\end{eqnarray*}
${\bf P}^\epsilon \mbox{-a.e.}$,
where $({\bf  M}_t(\epsilon, f))_{0 \le t \le T}$ is a continuous
$({\bf P}^\epsilon, ({\bf F}_t)_{0 \le t \le T})$-martingale
and $({\bf M}_t(\epsilon,$ $f)(r_{\scriptscriptstyle{T}}))_{0 \le t \le T}$
is a continuous $({\bf P}^\epsilon,
(r_{\scriptscriptstyle{T}}^{-1}({\bf F}_t))_{0 \le t \le T})$-martingale.
(We note that ${\bf P}^\epsilon \circ r_{\scriptscriptstyle{T}}^{-1}
= {\bf P}^\epsilon$ because $(T_{\epsilon,t})_{t \ge 0}$ is symmetric on
$L^2(\mu_\epsilon)$.) Moreover, by (\ref{eq50}) the bracket of ${\bf M}(\epsilon, f)$ is
given by
\begin{align*}
\langle {\bf M}(\epsilon, f) \rangle_t & = 2 \int_0^t
\epsilon^{d/2} \langle |\nabla f|_{{\mathbb R}^d}^2, {\bf X}(u) \rangle
+ \rho_{\phi}^{(1)}(\beta,1) \int_{{\mathbb R}^d}
|\nabla f(x)|_{{\mathbb R}^d}^2 \,dx \,du,
\end{align*}
as e.g.~directly follows from \cite{FOT94}, Theorem 5.2.3 and
Theorem 5.1.3(i). We note here that both theorems in \cite{FOT94}
are formulated and proved for locally compact separable metric spaces,
while ${\mathcal D}'$ is not of this type. However, both theorems carry over
to general state spaces by virtue of the local compactification and
regularization procedure developed in \cite{MaRo92}, Chap.~VI.2, which
is easy to see to be applicable in our case, see e.g.~\cite{MaRo92}, Chap.~VI,
Theorem 2.4, in regard to \cite{FOT94}, Theorem 5.1.3(i).
Hence by the Burkholder--Davies--Gundy inequalities and since
${\bf P}^\epsilon \circ r_{\scriptscriptstyle{T}}^{-1}
= {\bf P}^\epsilon$ we can find $C^{(3)} \in (0, \infty)$
such that for all $f \in {\mathcal D}, \, 0 < \epsilon \le 1, \, 0 \le s \le t \le
T$,
\begin{multline}\label{eq1605}
{\mathbb E}_{{\bf P}^\epsilon}[| \langle f, {\bf  X}(t) \rangle
- \langle f, {\bf X}(s) \rangle|^4]
\\  \le {\mathbb E}_{{\bf P}^\epsilon}[|{\bf  M}_t(\epsilon, f)
- {\bf M}_s(\epsilon, f)|^4 ]
+ {\mathbb E}_{{\bf P}^\epsilon}[
|{\bf  M}_{T-t}(\epsilon, f)(r_{\scriptscriptstyle{T}})
- {\bf M}_{T-s}(\epsilon, f)(r_{\scriptscriptstyle{T}})|^4 ]
\\  \le  C^{(3)} \, \Bigg{(} {\mathbb
E}_{{\bf P}^\epsilon}\Big{[} \Big{(}\int_s^{t} \Big{(}
\epsilon^{d/2} \langle |\nabla f|_{{\mathbb R}^d}^2,
{\bf X}(u) \rangle
+ \rho_{\phi}^{(1)}(\beta,1) \int_{{\mathbb R}^d}
|\nabla f(x)|_{{\mathbb R}^d}^2 \,dx \Big{)}
\,du \Big{)}^2 \Big{]} \\  +  {\mathbb E}_{{\bf
P}^\epsilon}\Big{[} \Big{(}\int_{T-t}^{T-s} \Big{(} \epsilon^{d/2}
\langle |\nabla f|_{{\mathbb R}^d}^2, {\bf X}(T-u)) \rangle
+ \rho_{\phi}^{(1)}(\beta,1) \int_{{\mathbb R}^d}
|\nabla f(x)|_{{\mathbb R}^d}^2 \,dx \Big{)}
\,du \Big{)}^2 \Big{]} \Bigg{)} \\
\le  4 \, C^{(3)} \,
(t-s)^{2} \Bigg{(} \epsilon^d \int_{\Gamma_{\epsilon}}
\langle |\nabla f|_{{\mathbb R}^d}^2, \omega \rangle^2 \,d\mu_{\epsilon}(\omega)
+ \rho_{\phi}^{(1)}(\beta,1)^2 \Big{(}\int_{{\mathbb R}^d}
|\nabla f(x)|_{{\mathbb R}^d}^2 \,dx \Big{)}^2 \Bigg{)} \\
\le  C^{(4)} \, (t-s)^{2} (\parallel |\nabla f|^2_{{\mathbb R}^d} \parallel^2_{d+1}
+ \parallel |\nabla f|_{{\mathbb R}^d} \parallel^4_0),
\end{multline}
where $C^{(4)} := 4 \, C^{(3)} \max( C^{(1)}), \rho^{(1)}_\phi(\beta,1)^{2}) $
and $C^{(1)}$ as in Theorem \ref{th9}{\rm (i)}.

Now we can use (\ref{eq1605}) to define $\langle f, {\bf  X}(t) \rangle
- \langle f, {\bf X}(s) \rangle$ for $f \in S({\mathbb R}^d)$ via an
approximation as an element in $L^4(\Omega, {\bf P}^\epsilon)$. Then, of course,
the estimate (\ref{eq1605}) is also true for $f \in S({\mathbb R}^d)$.

Let $m \in {\mathbb N}$ and let $(e_i)_{i \in {\mathbb N}}$ be the
sequence of Hermite functions, forming an orthonormal system in
${\mathcal H}_{m-2d}$. Then $(a_i^{m-d}e_i)_{i \in {\mathbb N}}$,
where $(a_i)_{i \in {\mathbb N}}$ are the eigenvalues of $A$
w.r.t.~the Hermite functions, forms an orthonormal system in ${\mathcal
H}_{-m}$. Since the mappings
$f \mapsto \parallel |\nabla f|^2_{{\mathbb R}^d} \parallel^2_{d+1}$ and
$f \mapsto \parallel |\nabla f|_{{\mathbb R}^d} \parallel^4_{0}$
are continuous on $S({\mathbb R}^d)$, we can choose $\alpha > 0$ and
$m \in {\mathbb N}$ large enough so that
\begin{eqnarray*}
\parallel |\nabla f|^2_{{\mathbb R}^d} \parallel^2_{d+1}
+ \parallel |\nabla f|_{{\mathbb R}^d} \parallel^4_0 \le \alpha
\parallel f \parallel^4_{m-2d}, \qquad \forall f \in S({\mathbb R}^d).
\end{eqnarray*}
In particular, for all $i \in {\mathbb N}$ we have
$\parallel |\nabla e_i|^2_{{\mathbb R}^d} \parallel^2_{d+1}
+ \parallel |\nabla e_i|_{{\mathbb R}^d} \parallel^4_0 \,\, \le \, \alpha.$
Hence, by the above we can estimate
\begin{multline}\label{eq95}
 \Big{(}{\mathbb E}_{{\bf P}^\epsilon}[\parallel {\bf X}(t)- {\bf X}(s)
\parallel^4_{-m}]\Big{)}^{1/2}
\!\!\!\!\!\! \le  \sum_{i = 0}^{\infty}  a_i^{-2d} \Big{(}{\mathbb E}_{{\bf P}^\epsilon} {[}
|\langle e_i, {\bf X}(t) \rangle
- \langle e_i, {\bf X}(s) \rangle|^4 {]} \Big{)}^{1/2} \!\!\!\!\!\!\le  C^{(5)} \, (t-s),
\end{multline}
where the constant $C^{(5)} :=  (\alpha C^{(4)})^{1/2}
\sum_{i = 0}^{\infty} a_i^{-2d}$ is finite, because
$A^{-d}$ is a Hilbert-Schmidt operator.
Since by Theorem \ref{th9}{\rm (iii)} $\mu_\epsilon \to \nu_\mu$ as $\epsilon \to 0$,
now the tightness of $({\bf P}^\epsilon)_{\epsilon > 0}$ on
$C([0, \infty), {\mathcal H}_{-m})$ follows by standard arguments.
\hfill $\blacksquare$

\subsection{Identification of the limit via the associated martingale problem}

In order to identify the limit by Theorem \ref{th2000} below it would
be sufficient to show that each accumulation point ${\bf P}$ of $({\bf
P}^\epsilon)_{\epsilon > 0}$ solves the martingale problem for $(-H,$
$D_0)$, where $D_0 := \{ G(\langle f, \cdot \rangle) \, | \,
G \in C^2_b({\mathbb R}), \, f \in S({\mathbb R}^d) \}$, with initial
distribution $\nu_\mu$, i.e.,
\begin{eqnarray*}
G(\langle f, {\bf X}(t) \rangle) - G(\langle f, {\bf X}(0) \rangle)
+ \int_0^t H G(\langle f, \cdot \rangle)({\bf X}(s)) \,ds, \qquad t \ge 0,
\end{eqnarray*}
is an ${\bf F}_t$-martingale under ${\bf P}$ and
${\bf P} \circ {\bf X}(0)^{-1} = \nu_\mu$.
One well-known way to establish this property is to prove
convergence of the generators $H_\epsilon$ to the generator
$H$ as $\epsilon \to 0$. Thus, first we study the difference
$\parallel (H - H_\epsilon)G(\langle f, \cdot \rangle)
\parallel_{L^2(\mu_\epsilon)}$
for $\epsilon \to 0$. In order see that $H G(\langle f, \cdot \rangle) \in
L^2(\mu_\epsilon)$ we use representation (\ref{eq5}).

\subsubsection{Non-convergence of generators}

Using (\ref{eq1147}) and (\ref{eq5}) again, by an approximation
argument it is easy to show that $\langle f, \cdot \rangle, \,
f \in {\mathcal D}$, is an
element of $D(H_\epsilon)$ and $D(H)$. As we shall prove now at least
on such functions the above convergence does not hold if we have
non-trivial interactions. For the proof of the following Theorem
we refer to Appendix \ref{a3}.

\begin{theorem}\label{th2001}
Let the potential $\phi$ be isotropic, i.e.,
$\phi(x) = V(r), \, r = |x|_{{\mathbb R}^d}$,
$x \in {\mathbb R}^d$. Furthermore,
let $x^k x^l \partial_i \partial_j \phi \in L^1({{\mathbb R}^d}, dx)$
and $x^i \partial_i \phi \in L^2({{\mathbb R}^d}, dx)$.
Additionally, let the assumptions required in Theorem \ref{pr200}
and Theorem \ref{dfdrr} below hold
and let $\mu$ be the Gibbs measure associated to $(\phi, \beta, 1)$
and the construction with empty boundary condition, where $\beta \in
[0, \beta_0]$ and $\beta_0 > 0$ is as in Theorem \ref{pr200}.
Then there exists a function $[0,\beta_0] \ni \beta
\mapsto R_\phi(\beta) \in {\mathbb R}_+$ such that
\begin{eqnarray*}
\lim_{\epsilon \to 0} \parallel (H - H_\epsilon)\langle f, \cdot
\rangle \parallel_{L^2(\mu_\epsilon)} = R_{\phi}(\beta)
\parallel \Delta f \parallel^2_{L^2(dx)}, \qquad \forall f \in {\mathcal D}.
\end{eqnarray*}
Furthermore, if $\mu \neq \pi_1$, then there exist $\beta_1(\phi)
\in (0, \beta_0]$ such that $R_\phi(\beta) > 0$
for all $\beta \in (0, \beta_1]$.
\end{theorem}

\begin{remark}{\rm Theorem \ref{th2001} states that for high
temperatures (small inverse temperature) and sufficiently smooth isotropic
potentials the generators do not converge in the $L^2$-sense.
It applies obviously to compactly supported potentials $\phi \in C_0^2({\mathbb
R}^d)$ and has been conjectured in \cite{Br80}, \cite{Ro81}, and
\cite{Sp86}.}
\end{remark}

\subsubsection{A conditional theorem on convergence in law}

In order to identify the limit the following weaker type
of convergence is sufficient.

\begin{conjecture}\label{le2000}
Let $(\phi, \beta, 1)$ fulfill conditions {\rm (SS)}, {\rm
(UI)}, {\rm (LR)}, {\rm (D)}, {\rm (LS)} and let $\mu$
be the corresponding Gibbs measure constructed with
empty boundary condition. Furthermore, for $G
\in C^2_b({\mathbb R})$, $f \in {\mathcal D}$, and $t, s \ge
0$, define
\begin{eqnarray*}
V_{\epsilon}(f,t,s) & := &  \int_t^{t+s}
G^{\prime}(\langle f, {\bf X}(u) \rangle)
(H - H_\epsilon)\langle f, \cdot \rangle({\bf X}(u)) \,du.
\end{eqnarray*}
Then
\begin{eqnarray*}
\lim_{\epsilon \to 0}{\mathbb
E}_{{\bf P}^\epsilon}[|V_{\epsilon}(f,t,s)|] = 0.
\end{eqnarray*}
\end{conjecture}

\begin{remark}\label{re0621}
{\rm Conjecture \ref{le2000} states that the scaled generators converge in time
average, whereas Theorem \ref{th2001} concerns convergence of the
scaled generators at an arbitrary fixed time.
Conjecture \ref{le2000} first has been formulated
in \cite{Ro81}.
In \cite{Sp86} the author describes a proof of Conjecture \ref{le2000}
for positive, smooth, compactly supported potentials and $d \le 3$, but
with $G(x) = x$, see \cite{Sp86}, assumption
(C), page 10. It is easy to show that, if Conjecture \ref{le2000} holds for
$G(x) = x$, then it also holds for all $G \in C^2_b({\mathbb R})$.}
\end{remark}

\begin{theorem}\label{th2000}
Let $(\phi, \beta, 1)$ fulfill conditions {\rm (SS)}, {\rm
(UI)}, {\rm (LR)}, {\rm (D)}, {\rm (LS)} and let $\mu$ be the
corresponding Gibbs measure constructed
with empty boundary condition. Assume
Conjecture {\rm \ref{le2000}}. Additionally, let ${\bf
P}$ be an accumulation point of $({\bf P}^\epsilon)_{\epsilon > 0}$
on $C([0, \infty), {\mathcal H}{-m})$ with $m \in {\mathbb N}$ as in Theorem
\ref{pr777}.
Then ${\bf P}$ solves the martingale problem for $(-H,$
$D_0)$ with initial distribution $\nu_\mu$, i.e.,
for all $G \in C^2_b({\mathbb R}), \, f \in S({\mathbb R}^d)$,
\begin{eqnarray}\label{eq29}
G(\langle f, {\bf X}(t)\rangle) - G(\langle f, {\bf X}(0) \rangle)
+ \int_0^t HG(\langle f, \cdot \rangle)({\bf X}(s)) \,ds, \qquad t \ge 0,
\end{eqnarray}
is an ${\bf F}_t$-martingale under ${\bf P}$ and
${\bf P} \circ {\bf X}(0)^{-1} = \nu_{\mu}$.
The measure ${\bf P}$ is uniquely determined by these properties,
in particular, all such ${\bf P}$ coincide. Hence
${\bf P}_{\epsilon} \to {\bf P}$ weakly as $\epsilon \to 0$.
\end{theorem}
{\bf Proof:}
Let $f \in {\mathcal D}, \, t, s \ge 0$, and define the following random variables on
$C([0, \infty)$, ${\mathcal H}_{-m})$:
\begin{align*}
U_{\epsilon}(f,t,s)  & :=   G(\langle f, {\bf X}(t) \rangle)
- G(\langle f, {\bf X}(s) \rangle) +
\int_t^{t+s} H_\epsilon G(\langle f, \cdot \rangle)({\bf X}(u))
\,du, \\ U(f,t,s) & :=  G(\langle f, {\bf X}(t) \rangle) -
G(\langle f, {\bf X}(s) \rangle) +
\int_t^{t+s} H G(\langle f, \cdot \rangle)({\bf X}(u))\,du, \\
S_{\epsilon}(f,t,s)  & := \epsilon^{d/2} \int_t^{t+s} G^{\prime
\prime}(\langle f, {\bf X}(u) \rangle) {\bf X}(|\nabla f|^2_{{\mathbb R}^d})(u)
\,du.
\end{align*}
Utilizing Theorem \ref{th9}{\rm (i)} it follows that
\begin{eqnarray}\label{eq43}
\lim_{\epsilon \to 0}{\mathbb
E}_{{\bf P}^\epsilon}[|S_{\epsilon}(f,t,s)|] = 0.
\end{eqnarray}
The trace filtration obtained by restricting $({\bf F}_t)_{t \ge
0}$ to $C([0, \infty), {\mathcal H}_{-m})$ coincides with the natural
filtration of $C([0, \infty), {\mathcal H}_{-m})$, which we also
denote by $({\bf F}_t)_{t \ge 0}$. Since
${\bf P}^\epsilon$ solves the martingale problem for
$(-H_{\epsilon} , D_0)$ w.r.t.~$({\bf F}_t)_{t \ge 0}$
we have for all ${\bf F}_t$-measurable bounded, continuous, $F_t: C([0,
\infty), {\mathcal H}_{-m}) \to {\mathbb R},$ and $\epsilon > 0$ that
${\mathbb E}_{{\bf P}^\epsilon}[F_t U_{\epsilon}(f,t,s)] = 0$.
Thus, together with Conjecture \ref{le2000} and (\ref{eq43}), it follows
that
\begin{eqnarray}\label{eq0520}
\lim_{\epsilon \to 0}{\mathbb
E}_{{\bf P}^\epsilon}[F_t U(f,t,s)]
= \lim_{\epsilon \to 0}{\mathbb
E}_{{\bf P}^\epsilon}[F_t (U_{\epsilon}(f,t,s) \!+\! V_{\epsilon}(f,t,s) \!+\!
S_{\epsilon}(f,t,s))] = 0.
\end{eqnarray}

Let ${\bf P}$ be an accumulation point of
$({\bf P}^\epsilon)_{\epsilon > 0}$ on $C([0, \infty), {\mathcal H}_{-m})$,
i.e., ${\bf P}^{\epsilon_n} \to {\bf P}$ weakly for some subsequence
$\epsilon_n \to 0$ for $n \to \infty$. Obviously, by Theorem \ref{th9}(iii)
we have ${\bf P} \circ {\bf X}(t)^{-1} = \nu_{\mu}$ for all $t \ge 0$, in
particular ${\bf P} \circ {\bf X}(0)^{-1} = \nu_{\mu}$. By (\ref{eq0520})
it remains to show that
\begin{eqnarray}\label{eq1520}
\lim_{n \to \infty}{\mathbb E}_{{\bf P}^{\epsilon_n}}[F_t U(f,t,s)]
= {\mathbb E}_{{\bf P}}[F_t U(f,t,s)].
\end{eqnarray}

Obviously, we only have to prove (\ref{eq1520}) with $U(f,t,s)$ replaced
by the last summand in its definition, because for the first two
summands convergence is clear. In order to do this we set
\begin{eqnarray*}
h := HG(\langle f, \cdot \rangle) =
- \rho_{\phi}^{(1)}(\beta,1)
G'' (\langle f,  \cdot \rangle) \parallel |\nabla f|_{{\mathbb R}^d}
\parallel_0^2
- \frac{\rho_{\phi}^{(1)}(\beta,1)}{\chi_{\phi}(\beta) }
G'(\langle f,  \cdot \rangle) \langle \Delta f,  \cdot \rangle.
\end{eqnarray*}
Then
\begin{multline*}
\Big{|}{\mathbb E}_{{\bf P}}\Big{[}F_t \int_t^{t+s}
H G(\langle f, \cdot \rangle)({\bf X}(u))\,du \Big{]} -
{\mathbb E}_{{\bf P}^{\epsilon_n}}\Big{[}F_t \int_t^{t+s}
H G(\langle f, \cdot \rangle)({\bf X}(u))\,du \Big{]} \Big{|} \\
\le \int_t^{t+s} {|}{\mathbb E}_{{\bf P}} {[}F_t
h({\bf X}(u)) {]} -
{\mathbb E}_{{\bf P}^{\epsilon_n}}{[}F_t
h({\bf X}(u)) {]} {|} \,du
\end{multline*}
and for $K_r := \{ \omega \in {\mathcal H}_{-m} \, | \parallel \omega
\parallel_{-m} \le r \}, \, r > 0,$ we have both for the positive and
negative parts $h^+, h^-$ of $h$ and $u \in [t, t+s]$, setting
$h_r^{\pm} := h^{\pm} \wedge \sup_{K_r}|h|,$
\begin{multline*}
{|}{\mathbb E}_{{\bf P}} {[}F_t
h^{\pm}({\bf X}(u)) {]} -
{\mathbb E}_{{\bf P}^{\epsilon_n}}{[}F_t
h^{\pm}({\bf X}(u)) {]} {|} \\
\le \Big{|} \int_{\{ {\bf X}(u) \in K_r \}}
|F_t| h_r^{\pm}({\bf X}(u)) \,d {\bf P}
- \int_{\{ {\bf X}(u) \in K_r \}}
\, |F_t| h_r^{\pm}({\bf X}(u)) \,d {\bf P}^{\epsilon_n} \Big{|} \\
+ \int_{ \{ {\bf X}(u) \in {\mathcal H}_{-m} \setminus K_r \}}
|F_t| |h|({\bf X}(u)) \,d {\bf P}
+ \int_{\{ {\bf X}(u) \in {\mathcal H}_{-m} \setminus K_r \}}
\! |F_t| |h|({\bf X}(u)) \,d {\bf P}^{\epsilon_n} \\
\le {|}{\mathbb E}_{{\bf P}} {[}|F_t|
h_r^{\pm}({\bf X}(u)) {]} -
{\mathbb E}_{{\bf P}^{\epsilon_n}}{[}|F_t|
h_r^{\pm}({\bf X}(u)) {]} {|} \\
+ 2 \int_{ \{ {\bf X}(u) \in {\mathcal H}_{-m} \setminus K_r \}}
\!\!\!\!\! |F_t| |h|({\bf X}(u)) \,d {\bf P}
+ 2 \int_{\{ {\bf X}(u) \in {\mathcal H}_{-m} \setminus K_r \}}
\!\!\!\!\! |F_t| |h|({\bf X}(u)) \,d {\bf P}^{\epsilon_n}
\end{multline*}
But for all $r > 0$
\begin{multline*}
\int_{\{ {\bf X}(u) \in {\mathcal H}_{-m} \setminus K_r \} }
|F_t| |h|({\bf X}(u)) \,d {\bf P}^{\epsilon_n}
\le \parallel F_t \parallel_{\infty}
\int_{{\mathcal H}_{-m} \setminus K_r } |h| \,d\mu_{\epsilon} \\
\le \frac{1}{r} \rho_{\phi}^{(1)}(\beta,1) \, C^{(1)}
\parallel F_t \parallel_{\infty}
\Big{(} \frac{\parallel G'' \parallel_{\infty}}{r} \parallel
|\nabla f|_{{\mathbb R}^d} \parallel^2_{0} +
\frac{\parallel G' \parallel_{\infty}}{\chi_{\phi}(\beta)}
\parallel
\Delta f \parallel_{m} \Big{)},
\end{multline*}
where we used $|\langle \Delta f, \omega \rangle| \le \|\Delta f\|_m
\|\omega\|_{-m}$ and $1 \le \|\omega\|_{-m}/r$ on
${\mathcal H}_{-m} \setminus K_r$. The constant
$C^{(1)}$ is as in Theorem \ref{th9}(i). Similarly,
\begin{multline*}
\int_{\{ {\bf X}(u) \in {\mathcal H}_{-m} \setminus K_r \} }
|F_t| |h|({\bf X}(u)) \,d {\bf P} \\
\le \rho_{\phi}^{(1)}(\beta,1)
\parallel F_t \parallel_{\infty}
\parallel G'' \parallel_{\infty} \parallel |\nabla f|_{{\mathbb R}^d}
\parallel^2_{0}  \frac{1}{r^2} \int_{{\mathcal H}_{-m}}
\parallel \omega \parallel^2_{-m} \,d\nu_{\mu}(\omega) \\
+ \frac{\rho_{\phi}^{(1)}(\beta,1)}{\chi_{\phi}(\beta) }
\parallel F_t \parallel_{\infty}  \parallel G' \parallel_{\infty}
\parallel \Delta f \parallel_{m}
\frac{1}{r}\int_{{\mathcal H}_{-m}}
\parallel \omega \parallel^2_{-m} \,d\nu_{\mu}(\omega),
\end{multline*}
and since the Gaussian measure $\nu_{\mu}$ has measure $1$ on
${\mathcal H}_{-m}$ there exists a constant $C^{(6)} \in (0,\infty)$
such that
$\int_{{\mathcal H}_{-m}}
\parallel \omega \parallel^2_{-m} \,d\nu_{\mu}(\omega) \le C^{(6)}$.
Hence by the weak convergence of ${\bf P}^{\epsilon_n} \to {\bf P}$
as $n \to \infty$ and Lebesgue's dominated convergence theorem
\begin{multline*}
\limsup_{n \to \infty}
\int_t^{t+s} |{\mathbb E}_{{\bf P}} [ F_t h^{\pm}({\bf X}(u)) ] -
{\mathbb E}_{{\bf P}^{\epsilon_n}}[F_t h^{\pm}({\bf X}(u))] | \,du \\
\le \frac{2 s}{r} \rho_{\phi}^{(1)}(\beta,1) \max\{C^{(1)}, C^{(6)}\}
\parallel F_t \parallel_{\infty}
\Big{(} \frac{\parallel G'' \parallel_{\infty}}{r} \parallel
|\nabla f|_{{\mathbb R}^d} \parallel^2_{0}
+ \frac{\parallel G' \parallel_{\infty}}{\chi_{\phi}(\beta) }
\parallel \Delta f \parallel_{m}\Big{)},
\end{multline*}
for all $r > 0$.
Letting $r \to \infty$ equality (\ref{eq1520}) follows and therefore
\begin{eqnarray}\label{eq1521}
{\mathbb E}_{{\bf P}}[F_t U(f,t,s)] = 0, \quad \forall f \in {\mathcal D}.
\end{eqnarray}
But by an approximation (\ref{eq1521}) is also true for all
$f \in S({\mathbb R}^d)$.

Now it remains to show that ${\bf P}$ is
uniquely determined by (\ref{eq29}). But this
follows by an easy generalization of Theorem 1.4 in \cite{HS78}.
All the assumptions required there are fulfilled in our situation
except for the assumption on the operator $B$. This operator $B$
in our case is $\sqrt{-\Delta}$, which is not bounded as required
in \cite{HS78}.
Analyzing the proof, however, one finds that continuity and
boundedness of the function
\begin{eqnarray*}
[0 , \infty) \ni t \mapsto \langle B \exp(t \Delta)f, B \exp(t \Delta)f
\rangle \in [0 , \infty)
\end{eqnarray*}
for a fixed $f \in S({\mathbb R}^d)$ is sufficient, which in our case
is obviously true.
\hfill $\blacksquare$

\begin{appendix}

\section{Inverse temperature derivative of correlation functions}

First, we have to define the finite volume
correlation functions
\begin{gather*}
\rho_{\phi, \Lambda}(\beta, z, \eta) :=  Z^{-1}_{\phi,
\Lambda}(\beta, z) \int_{\Gamma_{0, \Lambda}} \exp(- \beta
E^{\phi}_{\Lambda}(\eta \cup \xi)) \,d \lambda_z(\xi), \quad \beta
\ge 0, \, z > 0, \\ Z_{\phi, \Lambda}(\beta,z) := \int_{\Gamma_{0,
\Lambda}} \exp(- \beta E^{\phi}_{\Lambda}(\xi)) \,d \lambda_z(\xi),
\quad \eta \in \Gamma_{0, \Lambda}, \, \Lambda \in {\mathcal O}_c({\mathbb R}^{d}),
\end{gather*}
where we restricted the Lebesgue-Poisson measure to $\Gamma_{0,
\Lambda} := \bigsqcup_{n=0}^\infty \Gamma_{0, \Lambda}^{(n)}$, see
Section \ref{ss23}.

The proof of the following lemma is an easy generalization of
Theorem 3.3.18 in \cite{K99}.

\begin{lemma}\label{le0521}
Let $(\phi, \beta_0, z)$ satisfy conditions
{\rm (S)} and {\rm (UI)}. Furthermore, let $\phi$ fulfill the condition
\begin{eqnarray}\label{eq0523}
0 < \int_{{\mathbb R}^d \setminus \Lambda_0}
(\exp(\beta_0 |\phi(x)|) - 1) \,dx < \infty
\end{eqnarray}
for some $\Lambda_0 \in {\mathcal O}_c({\mathbb R}^d)$.
Then
\begin{eqnarray}\label{eq29a}
\lim_{\Lambda \nearrow {\mathbb R}^d} \rho_{\phi,
\Lambda}^{(n)}(\beta, z, x_1,\dots,x_n) =
\rho_\phi^{(n)}(\beta, z, x_1,\dots,x_n)
\end{eqnarray}
for all $z > 0$ and uniformly in
$\beta, x_1, \ldots x_n$ on any set
$[0, \beta_0] \times (\Lambda')^n$, where
$\Lambda' \in {\mathcal O}_c({\mathbb R}^d)$.
\end{lemma}

\begin{remark}
{\rm Condition (\ref{eq0523}) is obviously fulfilled for smooth,
compactly supported potentials $\phi$. Or, if
$\phi \in L^1({\mathbb R}^d \setminus \Lambda_0)$ and
bounded on ${\mathbb R}^d \setminus \Lambda_0$
for some $\Lambda_0 \, \in {\mathcal O}_c({\mathbb R}^d)$,
and not $dx$-a.e.~zero on ${\mathbb R}^d
\setminus \Lambda_0$, then condition (\ref{eq0523}) is also
fulfilled.}
\end{remark}

Via a recursion formula one can transform the correlation
functions $\rho^{(n)}_{\phi, \Lambda}$ into the so-called {\it Ursell
functions} $u^{(n)}_{\phi, \Lambda}$ and vice versa, see e.g. \cite{MM91},
\cite{Rue69}. Their relation is given by
\begin{eqnarray}\label{eq22a}
\rho_{\phi, \Lambda}(\beta, z, \eta) = \sum_{\stackunder{\eta_k \cap \eta_l
= \emptyset, k \neq l, j \in {\mathbb N}}{\eta_1 \cup \ldots \cup
\eta_j = \eta}} u_{\phi, \Lambda}(\beta, z, \eta_1) \dotsm
u_{\phi, \Lambda}(\beta, z, \eta_j), \qquad \eta \in \Gamma_0,
\end{eqnarray}
where $u_{\phi, \Lambda}^{(n)}$ is related to $u_{\phi, \Lambda}$ analogously
to (\ref{eq367}).
Correspondingly, $u_{\phi}^{(n)}$ and $u_{\phi}$ are defined with
$\rho_{\phi}$ replacing $\rho_{\phi, \Lambda}$.
Due to the translation invariance of the correlation functions,
Ursell functions are also translation invariant. Furthermore,
by an easy generalization of Theorem 4.5 in \cite{Br80}, see also
\cite{Rue69}, Chapter 4, we obtain the following integrability property.

\begin{proposition}\label{pr1969}
Let $(\phi, \beta, z)$ satisfy conditions {\rm (S)} and {\rm (UI)}.
Then for each $n \ge 1$, there exists a non-negative measurable function
$U^{(n+1)}_{\phi, \beta, z}: {\mathbb R}^{d \times n} \to {\mathbb R}_+$, such that
\begin{align*}
|u^{(n+1)}_{\phi, \Lambda}(\beta, z, \cdot, 0)| \le U^{(n+1)}_{\phi, \beta, z},
\quad \forall \Lambda \in {\mathcal O}_c({\mathbb R}^d),
\end{align*}
and
\begin{multline*}
\int_{{\mathbb R}^{d \times n}} |U_{\phi, \beta, z}^{(n+1)}(x_1, \ldots, x_n)
f(x_1, \ldots, x_n)|\,dx_1 \ldots \,dx_n \\
\le \exp(2n \beta B(\phi)) \Big{(} \sum_{m=0}^{\infty} \frac{1}{m!}
(n+m+1)^{n+m-1} C(\beta \phi, z)^m \Big{)} \\
\times \sup_{x_{n} \in {\mathbb R}^d}
\int_{{\mathbb R}^d} |\exp(-\beta \phi(x_{n}-y_{n}))-1|
\sup_{x_{n-1} \in {\mathbb R}^d}
\int_{{\mathbb R}^d} |\exp(-\beta \phi(x_{n-1}-y_{n-1}))-1| \\ \dotsm
\sup_{x_{1} \in {\mathbb R}^d} \int_{{\mathbb R}^d} |\exp(-\beta \phi(x_1-y_1))-1|
|f(y_1,\dots,y_n)| \,dy_1 \ldots \,dy_n,
\end{multline*}
for all measurable functions $f: {\mathbb R}^{d \times n} \to {\mathbb R}$.
\end{proposition}

\begin{theorem}\label{pr200}
Let $(\phi, \beta_0, z)$ satisfy conditions {\rm (S)}, {\rm (UI)}, and let
either $\phi \equiv 0$ or
$\phi \in L^1({\mathbb R}^d, dx) \cap L^{2}({\mathbb R}^d, dx)$
and condition (\ref{eq0523}) hold.
Then $\rho_{\phi} \in C^{2}([0, \beta_0])$ and for
$\lambda_z$-a.e.~$\eta \in \Gamma_0$ we have
\begin{multline}\label{eq22c}
\frac{\partial \rho_{\phi}}{\partial \beta}(\beta, z, \eta) = -
E^{\phi}(\eta) \rho_{\phi}(\beta, z, \eta) - \int_{{\mathbb R}^d}
W^{\phi}(\eta \mid x) \rho_{\phi}(\beta, z, \eta \cup \{ x \}) \, z \,dx \\
- \frac{1}{2} \int_{{\mathbb R}^{d \times 2}} \phi(x-y) \Big{(}
\rho_{\phi}(\beta, z, \eta \cup \{ x, y \}) -  \rho_{\phi}(\beta,
z, \eta)\rho^{(2)}_{\phi}(\beta, z, x, y) \Big{)} \, z^2 \,dx \,dy,
\end{multline}
where $E^{\phi}(\eta) := \lim_{\Lambda
\nearrow {\mathbb R}^d} E_{\Lambda}^{\phi}(\eta)$.
\end{theorem}
{\bf Proof:} First, we note that the expression on the r.h.s.~of
(\ref{eq22c}) is well-defined and finite. Indeed, since $\phi \in
L^1({\mathbb R}^d)$ and the correlation functions are bounded, see
(\ref{eq18}), the first integral in this expression is finite.
Using (\ref{eq22a}) and Proposition \ref{pr1969}, one finds that the
second integral is also finite.

Analyzing the properties of the Lebesgue-Poisson measure we find for
$\eta \in \Gamma_{0, \Lambda}$:
\begin{multline*}
\frac{\partial \rho_{\phi, \Lambda}}{\partial \beta}(\beta, z,
\eta) = - E_{\Lambda}^{\phi}(\eta) \rho_{\phi, \Lambda}(\beta, z,
\eta) - \int_{\Lambda} W^{\phi}(\eta \mid x) \rho_{\phi,
\Lambda}(\beta, z, \eta \cup \{ x \}) \, z \,dx \\ -  \frac{1}{2}
\int_{\Lambda^2} \phi(x-y) \Big{(} \rho_{\phi, \Lambda}(\beta, z,
\eta \cup \{ x, y \}) - \rho_{\phi, \Lambda}(\beta, z,
\eta)\rho^{(2)}_{\phi, \Lambda}(\beta, z, x,y) \Big{)} \, z^2 \,dx \,dy.
\end{multline*}

Using (\ref{eq29a}), Proposition \ref{pr1969}, and that
bound (\ref{eq18}) also holds for finite volume correlation functions,
uniformly in $\Lambda \subset {\mathbb R}^d$, twice applying
the dominated convergence theorem shows
\begin{multline*}
\lim_{\Lambda \nearrow {\mathbb R}^d} \!\Big{(} \frac{\partial
\rho_{\phi, \Lambda}}{\partial \beta}\!\Big{)}(\beta, z, \eta)  = \! -
E^{\phi}(\eta) \rho_{\phi}(\beta, z, \eta) - \!\int_{{\mathbb R}^d}
\!\!\!\! W^{\phi}(\eta \mid x) \rho_{\phi}(\beta, z, \eta \cup \{ x \})
z \,dx \\ - \frac{1}{2} \int_{{\mathbb R}^{d \times 2}} \phi(x-y) \Big{(}
\rho_{\phi}(\beta, z, \eta \cup \{ x, y \}) -  \rho_{\phi}(\beta,
z, \eta)\rho^{(2)}_{\phi}(\beta, z, x,y) \Big{)} \, z^2 \,dx \,dy,
\end{multline*}
for $\lambda_z$-a.e.~$\eta \in {\Gamma_0}$.
It remains to show that derivative and the infinite volume limit can
be interchanged. We evidently have to show this only for potentials
which are not identically equal to zero. By using Lemma \ref{le0521} and Proposition
\ref{pr1969}, we
see that for $z > 0, \, \eta \in \Gamma_0$, fixed the function
$\frac{\partial \rho_{\phi, \Lambda}}{\partial \beta}(\beta, z,
\eta)$ converges uniformly on $[0, \beta_0]$ as
$\Lambda \nearrow {\mathbb R}^d$ and
\begin{eqnarray*}
\frac{\partial \rho_{\phi}}{\partial \beta}(\beta, z, \eta)
= \lim_{\Lambda \nearrow {\mathbb R}^d}
\frac{\partial \rho_{\phi, \Lambda}}{\partial \beta}(\beta, z, \eta).
\end{eqnarray*}

The second order derivative can be derived analogously.
The only difference is that in the second order derivative the
potential $\phi$ appears in its second power. Hence, for
$\phi \in L^1({\mathbb R}^d, dx) \cap L^{2}({\mathbb R}^d, dx)$
we obtain that $\rho_{\phi} \in C^{2}([0, \beta_0])$.
A more detailed proof can be found in \cite{GKLR00}
\hfill$\blacksquare$

\section{Coercivity identity for Gibbs measures}

Here we derive an analog of the usual coercivity identity on
$L^2 ({\mathbb R}^d, g\,dx)$ for $L^2(\Gamma,\mu)$,
where $\mu$ is a Ruelle measure on $\Gamma$, whose potential
satisfies some weak additional conditions.

First we have to develop a little further the analysis and geometry
as in Section \ref{s44}. For each $\gamma \in \Gamma$, consider the triple
\begin{eqnarray}\label{eq1405}
T_{\gamma ,\,\infty }(\Gamma) \supset T_\gamma (\Gamma)\supset
T_{\gamma ,0}(\Gamma).
\end{eqnarray}
Here, $T_{\gamma ,0}(\Gamma )$ consists of all finite sequences
from $T_\gamma (\Gamma) $, and $T_{\gamma ,\,\infty }(\Gamma) :=
\left( T_{\gamma ,0}(\Gamma) \right)^{\prime }$ is the dual space
consisting of all sequences $V(\gamma )=(V(\gamma ,x))_{x\in
\gamma }$, where $V(\gamma ,x)\in T_x({\mathbb R}^d)$. The pairing
between any $V(\gamma )\in T_{\gamma ,\,\infty }(\Gamma) $ and
$v(\gamma )\in T_{\gamma ,0}(\Gamma) $ with respect to the zero
space $T_\gamma (\Gamma) $ is given by
\begin{eqnarray*}
( V(\gamma ),v(\gamma )) _{T_{\gamma}(\Gamma)} := \sum_{x\in
\gamma }( V(\gamma ,x),v(\gamma ,x)) _{T_x({\mathbb R}^d)}.
\end{eqnarray*}
This series is, in fact, finite.

For $\gamma\in\Gamma$, we define
$B_\mu(\gamma)=(B_\mu(\gamma,x))_{x\in\gamma}\in
T_{\gamma,\,\infty}(\Gamma)$ by
\begin{eqnarray}\label{453453}
B_\mu(\gamma,x) := -\beta \sum_{y\in\gamma\setminus\{x\}}\nabla
\phi(x-y), \qquad x \in \gamma.
\end{eqnarray}
As follows from the proof of Lemma 4.1 in \cite{AKR98b}, for
$\mu$-a.e.~$\gamma\in\Gamma$ the series on the right hand side of
(\ref{453453}) converges absolutely in ${\mathbb R}^d$, provided
$(\phi, \beta, z)$ satisfies (SS), (UI), (LR) and (D),
and that $\mu$ is the corresponding Gibbs measure constructed with
empty boundary condition.
Observe, that

\begin{gather}\label{eq4378}
H_\mu^\Gamma F(\gamma) := -\Delta^\Gamma F(\gamma)
-(B_\mu(\gamma),\nabla^\Gamma F(\gamma))_{T_{\gamma}(\Gamma)},  \\
\Delta^\Gamma F(\gamma) := \sum_{x\in\gamma}\Delta_xF(\gamma),
\quad \Delta_x F(\gamma) := \Delta_y F_x(x,y)\big|_{y=x},
\nonumber
\end{gather}
where $H_\mu^\Gamma$ is the generator as in \eqref{eq37} and
$F\in {\mathcal F}C_{b}^\infty ({\mathcal D}, \Gamma )$ as in (\ref{eq54}).
Of course, $\Delta$ acting on differentiable functions
defined on ${\mathbb R}^d$ is denoting the Laplacian
on ${\mathbb R}^d$.
We call $B_\mu$ the logarithmic derivative of
the measure $\mu$.

Let $A(\gamma) \in (T_{\gamma,\,\infty}(\Gamma))^{\otimes 2}$,
cf.~(\ref{eq1405}), so that
$A(\gamma)=(A(\gamma,x,y))_{x,y\in\gamma}$, where
$A(\gamma,x,y)$ $\in T_y({\mathbb R}^d)\otimes T_x({\mathbb R}^d)$.
We realize $A(\gamma)$ as a linear operator acting from
$T_{\gamma,0}(\Gamma)$ into $T_{\gamma,\,\infty}(\Gamma)$ setting
\begin{gather*}
T_{\gamma,0}(\Gamma) \ni V(\gamma) \mapsto
A(\gamma)V(\gamma) \, :=
\Big(\sum_{x\in\gamma}(A(\gamma,x,y),V(\gamma,x))_{T_x({\mathbb
R}^d)}\Big)_{y\in\gamma}\!\!\!\!\!\in T_{\gamma,\,\infty}(\Gamma).
\end{gather*}
Evidently, if $A(\gamma)\in (T_{\gamma,0}(\Gamma))^{\otimes 2}$,
then $A(\gamma)$ defines a linear continuous operator in $T_\gamma
(\Gamma)$. We denote by $A(\gamma)^*$ its adjoint operator.

For a vector field $\Gamma\ni\gamma\mapsto W(\gamma)\in
T_{\gamma,\,\infty}(\Gamma)$, we define its derivative
$\nabla^\Gamma W(\gamma)$ as a mapping
\begin{gather*}
\Gamma\ni\gamma \mapsto
\nabla^\Gamma W(\gamma)=(\nabla^\Gamma
W(\gamma,x,y))_{x,y\in\gamma}\in
(T_{\gamma,\,\infty}(\Gamma))^{\otimes2}
\end{gather*}
 such that
\begin{gather*}
\nabla^\Gamma W(\gamma,x,y){:=}\nabla_y W (\gamma,x)=\begin{cases}
\nabla_z W (\gamma-\varepsilon_y+\varepsilon_z,x)\big|_{z=y},&
\text{if }x\ne y,\\ \nabla_z
W(\gamma-\varepsilon_y+\varepsilon_z,z)\big|_{z=y},& \text{if
}x=y,\end{cases}
\end{gather*}
 if all  derivatives $\nabla_y W(\gamma,x)$,
$x,y\in\gamma$, exist. For a function $F\colon\Gamma\to{\mathbb R}^d$, we
denote $F''{:=}\nabla^\Gamma\nabla^\Gamma F$, if it exists.

\begin{theorem}[coercivity identity]\label{dfdrr}
Let the potential $\phi$ satisfy {\rm (SS)}, $\!${\rm (I)}, {\rm (LR)}, and the three
following conditions:

\noindent
{\rm (i)} $\phi\in C^2({\mathbb R}^d\setminus\{0\})$,
$e^{-\phi}$ is continuous on ${\mathbb R}^d$, and
$e^{-\phi}\nabla\phi$ can be extended to a continuous,
vector-valued function on ${\mathbb R}^d$;

\noindent {\rm (ii)} for each $\gamma\in S_\infty$, the three series
$\sum_{x\in\gamma}\phi(\cdot-x)$,
$\sum_{x\in\gamma}\nabla\phi(\cdot-x)$, and
$\sum_{x\in\gamma}\nabla^2\phi(\cdot-x)$ converge locally uniformly on
$X\setminus \gamma$;

\noindent{\rm (iii)}
we have
\begin{gather*}
\nabla\phi \in L^1({\mathbb
R}^d,\exp(-\phi(x))\, dx)\cap L^2({\mathbb
R}^d,\exp(-\phi(x))\, dx),\\  \nabla^2\phi\in
L^1({\mathbb R}^d,\exp(-\phi(x))\, dx).
\end{gather*}
Furthermore, let $\mu$ be the Gibbs measure corresponding to
$(\phi, \beta, z)$ and the construction with empty boundary
condition. Then, for any
$F\in {\mathcal F}C_b^{\infty}({\mathcal D}, \Gamma)$:
\begin{multline}\label{ewweqwq}
\| H_\mu^\Gamma
F\|^2_{L^2(\mu)}=\int_{\Gamma}\operatorname{Tr}_{T_\gamma(\Gamma)}F''(\gamma)F''(\gamma)^*\,d\mu(\gamma)
-\int_\Gamma (\nabla^\Gamma
F(\gamma),\nabla^\Gamma B_\mu(\gamma)\nabla^\Gamma
F(\gamma))_{T_\gamma(\Gamma)}\,d\mu(\gamma)\\
=\int_{\Gamma}\operatorname{Tr}_{T_\gamma(\Gamma)}F''(\gamma)F''(\gamma)^*\,d\mu
(\gamma)+ \beta \int_\Gamma \sum_{\{x,y\}\subset\gamma
}\big((\nabla^\Gamma F(\gamma,x)-\nabla^\Gamma
F(\gamma,y)),\\
\nabla^2\phi(x-y) (\nabla^\Gamma
F(\gamma,x)-\nabla^\Gamma F(\gamma,y))\big)_{{\mathbb
R}^d}\,d\mu(\gamma).
\end{multline}
\end{theorem}

\begin{remark}{\rm As easily seen, conditions (i)--(iii) of the
above theorem imply (D) and (LS).}
\end{remark}

\begin{remark}
{\rm As will be seen from the proof of
Theorem \ref{dfdrr}, the coercivity identity \eqref{ewweqwq}
holds for each monomial $F=\langle f,\cdot \rangle^n$, where
$f\in{\mathcal D}$ and $n\in{\mathbb N}$.}
\end{remark}
{\bf Proof:} Let $G:\Gamma\times{\mathbb R}^d\to{\mathbb R}_+$ be
measurable, then, by \cite{NZ79}, we have due to condition {\rm (ii)}:
\begin{eqnarray}\label{frerws}
\int_\Gamma \sum_{x\in\gamma}
G(\gamma,x)\,d\mu(\gamma)=\!\!\int_{\Gamma}\!\int_{{\mathbb
R ^d}}\!\!z \,\exp\Big(\!\!\!-\!\!\beta\!\sum_{y\in\gamma} \!\! \phi(x-y)\Big)
G(\gamma+\varepsilon_x,x) \,dx \,d\mu(\gamma).
\end{eqnarray}

Let $F\in {\mathcal F}C_b^{\infty}({\mathcal D}, \Gamma)$. By
\eqref{eq4378} and (\ref{frerws}), we get
\begin{multline} \|H_\mu^\Gamma F\|^2_{L^2(\mu)}=\int_\Gamma \sum_{x\in\gamma}
\big( \Delta_x F(\gamma)+(B_\mu(\gamma,x),\nabla^\Gamma
F(\gamma,x))_{T_x({\mathbb R}^d)} \big)^2\,d\mu(\gamma) \\
\text{}+\int_{\Gamma}\sum_{x,y\in\gamma,\, x\ne y}\big( \Delta_x
F(\gamma)+(B_\mu(\gamma,x),\nabla^\Gamma
F(\gamma,x))_{T_x({\mathbb R}^d)} \big) \\ \times \big(
\Delta_y F(\gamma)+(B_\mu(\gamma,y),\nabla^\Gamma
F(\gamma,y))_{T_y({\mathbb R}^d)} \big)\,d\mu(\gamma)\\
=\int_\Gamma \int_{{\mathbb R}^d} z \,
\exp\Big(-\beta\sum_{y\in\gamma}\phi(x-y)\Big)\\ \times
\big(\Delta_x
F(\gamma+\varepsilon_x)+(B_\mu(\gamma+\varepsilon_x,x),\nabla_x
F(\gamma+\varepsilon_x))_{T_x({\mathbb R}^d)}\big)^2
\,dx \,d\mu(\gamma) \\
+ \int_\Gamma \int_{{\mathbb R}^d} \int_{{\mathbb R}^d} z^2 \,
\exp\Big(-\beta \Big( \sum_{y_1\in\gamma}\phi(x_1-y_1)
+ \sum_{y_2\in\gamma\cup\{x_1\}}\phi(x_2-y_2) \Big)\Big)\\
\times \big(\Delta_{x_1}
F(\gamma+\varepsilon_{x_1}\!+\varepsilon_{x_2})\!+\!(B_\mu(\gamma+
\varepsilon_{x_1}+\varepsilon_{x_2},x_1),\!\nabla_{x_1}
F(\gamma+\varepsilon_{x_1}+\varepsilon_{x_2}))_{T_{x_1}({\mathbb
R}^d)} \big)\\ \times \big(\Delta_{x_2}
F(\gamma+\varepsilon_{x_1}\!+\varepsilon_{x_2})\!+\!(B_\mu(\gamma+
\varepsilon_{x_1}+\varepsilon_{x_2},x_2), \\
\nabla_{x_2}
F(\gamma+\varepsilon_{x_1}+\varepsilon_{x_2}))_{T_{x_2}({\mathbb
R}^d)} \big) \,dx_1 \,dx_2 \,d\mu(\gamma). \label{ysawa}
\end{multline}

By  conditions {\rm (i)} and {\rm (ii)}, we conclude that, for each fixed
$\gamma\in S_\infty$, the function
\begin{eqnarray*}
g_\gamma(x) :=
\exp\Big(-\beta\sum_{y\in\gamma}\phi(x-y)\Big)
\end{eqnarray*}
is continuous
on ${\mathbb R}^d$, two times continuously differentiable on
${\mathbb R}^d\setminus\gamma$, and $\nabla g_\gamma$ extends to a
continuous function on ${\mathbb R}^d$. Moreover, by \eqref{453453},
$B_\mu(\gamma+\varepsilon_x,x)$ is the logarithmic derivative of the
measure $\nu_\gamma := g_\gamma\,dx$. Finally, it is easy to see from
{\rm (i)}--{\rm (iii)} that the function
\begin{eqnarray*}
g_\gamma(x)\big(\log
g_\gamma(x)\big)''=\exp\Big(-\beta\sum_{y\in\gamma}\phi(x-y)\Big)
\beta\sum_{y\in\gamma}\nabla^2\phi(x-y)
\end{eqnarray*}
 belongs to $L^1_{\rm loc}({\mathbb R}^d)$. Thus, the usual coercivity identity on
the space of square-integrable functions $L^2({\mathbb
R}^d,d\nu_\gamma)$ implies that
\begin{multline}
\int_{{\mathbb R }^d} \exp\Big(-\beta\sum_{y\in\gamma}
\phi(x-y)\Big)
\big(\Delta_x
F(\gamma+\varepsilon_x)+(B_\mu(\gamma+\varepsilon_x,x),\nabla_x
F(\gamma+\varepsilon_x))_{T_x({\mathbb R}^d)}\big)^2 \,dx  \\
=\int_{{\mathbb R }^d}
\exp\Big(-\beta\sum_{y\in\gamma}\phi(x-y)\Big)\big(\operatorname{Tr}_{T_x({\mathbb
R }^d)}\nabla_x\nabla_x F(\gamma+\varepsilon_x)(\nabla_x\nabla_x
F(\gamma+\varepsilon_x))^* \\ \text{}- (\nabla_x
F(\gamma+\varepsilon_x),\nabla_x
B_\mu(\gamma+\varepsilon_x,x)\nabla_x F(\gamma+\varepsilon_x)
)_{T_x({\mathbb R}^d)}\big) \,dx.
\end{multline}
Absolute analogously, a slight modification of the proof of the
coercivity identity on ${\mathbb R}^d$ implies that
\begin{multline}\label{resji}
\int_{{\mathbb R}^d}\int_{{\mathbb R}^d}\exp\Big(
-\beta\sum_{y_1\in\gamma}
\phi(x_1-y_1)-\beta\sum_{y_2\in\gamma}\phi(x_2-y_2)-\beta\phi(x_1-x_2)\Big)
\big(\Delta_{x_1}
F(\gamma+\varepsilon_{x_1}+\varepsilon_{x_2}) \\
+(B_\mu(\gamma+\varepsilon_{x_1}+\varepsilon_{x_2},x_1),\nabla_{x_1}
F(\gamma+\varepsilon_{x_1}+\varepsilon_{x_2}))_{T_{x_1}({\mathbb
R}^d)} \big) \big(\Delta_{x_2}
F(\gamma+\varepsilon_{x_1}+\varepsilon_{x_2}) \\
+(B_\mu(\gamma+\varepsilon_{x_1}+\varepsilon_{x_2},x_2),\nabla_{x_2}
F(\gamma+\varepsilon_{x_1}+\varepsilon_{x_2}))_{T_{x_2}({\mathbb
R}^d)} \big) \,dx_1 \,dx_2 \\ =\int_{{\mathbb
R}^d}\int_{{\mathbb R}^d} \,\exp\Big(
-\beta\sum_{y_1\in\gamma}
\phi(x_1-y_1)-\beta\sum_{y_2\in\gamma}\phi(x_2-y_2)-\beta\phi(x_1-x_2)\Big)
\\ \big(\|\nabla_{x_2}\nabla_{x_1} F(\gamma+ \varepsilon_{x_1}+
\varepsilon_{x_2})\|_{T_{x_2}({\mathbb R}^d)\otimes
T_{x_1}({\mathbb R}^d)}^2 - (\nabla_{x_2}
F(\gamma+\varepsilon_{x_1}+\varepsilon_{x_2}), \\ \nabla_{x_2} B_\mu
(\gamma+\varepsilon_{x_1}+\varepsilon_{x_2},x_1)\nabla_{x_1}
F(\gamma+\varepsilon_{x_1}+\varepsilon_{x_2}))_{T_{x_2}({\mathbb
R}^d)}\big) \,dx_1 \,dx_2.
\end{multline}

Next, by \eqref{453453}, {\rm (i)}, and {\rm (ii)}, we get for any $\gamma\in
S_\infty$:
\begin{eqnarray}\label{gfdtrstrs} \nabla_y B_\mu
(\gamma,x)=\begin{cases}-\beta\sum\limits_{z\in\gamma\setminus\{x\}}
\nabla^2\phi(x-z),&\text{if }x=y,\\
\beta\nabla^2\phi(x-y),&\text{otherwise.}\end{cases}
\end{eqnarray}
By \eqref{frerws}, \eqref{gfdtrstrs}, condition {\rm (iii)}, and estimate
(4.29) in \cite{AKR98b}, we have for any $\Lambda\in{\mathcal
O}_c({\mathbb R}^d)$:
\begin{gather} \int_\Gamma \int_\Lambda z \, \exp\Big(
\!-\!\beta\sum_{y\in\gamma}\phi(x-y)\Big) \big(1+\|\nabla_x
B_\mu(\gamma+\varepsilon_x,x)\|_{T_x({\mathbb R}^d)\otimes T_x({\mathbb R}^d)}\big)
\notag \,dx \,d\mu(\gamma)  \\ =\int_\Gamma \sum_{x\in\gamma_\Lambda}
\big(1+\|\nabla_x B_\mu(\gamma,x)\|_{T_x({\mathbb R}^d)\otimes
T_x({\mathbb R}^d)}\big) \,d\mu(\gamma) \notag \\ \le
\int_\Gamma \sum_{x\in\gamma_\Lambda}
\Big(1+\beta\sum_{y\in\gamma\setminus\{x\}}\|
\nabla^2\phi(x-y)\|_{{\mathbb R}^d\otimes{\mathbb R}^d}\Big) \,d\mu(\gamma)
\notag\\ =\int_\Lambda \rho_\mu^{(1)}(x) \,dx +\int_\Lambda
\int_{{\mathbb R}^d} \rho_\mu^{(2)}(x,y) \, \beta \,
\, \|\nabla^2\phi(x-y)\|_{{\mathbb R}^d\otimes {\mathbb R}^d}
\,dy \,dx \notag\\
\le \int_\Lambda \rho_\mu^{(1)}(x) \,dx
+ C^{(8)} \!\! \int_\Lambda \int_{{\mathbb R}^d} \beta
\, \|\nabla^2\phi(x-y)\|_{{\mathbb R}^d\otimes{\mathbb R}^d}
\,e^{-\beta\phi(x-y)} \,dy \,dx <\infty,\label{dewse}
\end{gather}
where $C^{(8)} \in (0, \infty)$ is a constant,
and analogously
\begin{gather} \int_\Gamma \int_\Lambda \int_\Lambda z^2
\exp\Big( -\beta\sum_{y_1\in\gamma}
\phi(x_1-y_1)-\beta\sum_{y_2\in\gamma}\phi(x_2-y_2)-\beta\phi(x_1-x_2)\Big)\notag\\
\times \big( 1+\|\nabla_{x_2}
B_\mu(\gamma+\varepsilon_{x_1}+\varepsilon_{x_2},x_1)\|_{T_{x_2}({\mathbb R}^d)\otimes
T_{x_1}({\mathbb R}^d)}\big) \,dx_1 \,dx_2 \,d\mu(\gamma) \notag\\ =
\int_\Gamma \sum_{x,y\in\gamma_\Lambda,\, x\ne
y}\big(1+\|\nabla_y B_\mu(\gamma,x)\|_{T_y({\mathbb R}^d)\otimes
T_x({\mathbb R}^d)}\big) \,d\mu(\gamma) \notag \\
= \int_\Gamma \sum_{x,y\in\gamma_\Lambda,\, x\ne y}
\big(1+\beta \, \|\nabla^2\phi(x-y)\|_{{\mathbb R}^d\otimes{\mathbb
R}^d} \big) \,d\mu(\gamma) <\infty. \label{wqwqw}
\end{gather}
Now, by \eqref{frerws}--\eqref{resji},
\eqref{dewse}, and \eqref{wqwqw},
\begin{multline}\|H_\mu^\Gamma F\|^2_{L^2(\mu)}=\int_\Gamma
\sum_{x\in\gamma}\big(
\operatorname{Tr}_{T_x({\mathbb R}^d)}\nabla_x\nabla_x F(\gamma)
(\nabla_x\nabla_x F(\gamma))^* \\
\vphantom{\sum_{x\in\gamma}}\text{}-(\nabla_x
F(\gamma),\nabla_x B_\mu(\gamma,x)\nabla_x
F(\gamma))_{T_x({\mathbb R}^d)}\big) +
\!\!\! \sum_{x,y\in\gamma,\,x\ne y} \!\!\!\Big(\|\nabla_y\nabla_x
F(\gamma)\|^2_{T_y({\mathbb R}^d)\otimes T_x({\mathbb R}^d)} \\
-(\nabla_y
F(\gamma),\nabla_y B_\mu(\gamma,x)\nabla_x
F(\gamma))_{T_y({\mathbb R}^d)}\Big) \,d\mu(\gamma)\\
=\!\int_{\Gamma}\operatorname{Tr}_{T_\gamma(\Gamma)}F''(\gamma)F''(\gamma)^*\,d\mu
(\gamma)-\!\int_\Gamma (\nabla^\Gamma F(\gamma),\nabla^\Gamma
B_\mu(\gamma)\nabla^\Gamma
F(\gamma))_{T_\gamma(\Gamma)}\,d\mu(\gamma). \label{wewepioi}
\end{multline}
Finally, from \eqref{gfdtrstrs} and \eqref{wewepioi} we get the
second equality in \eqref{ewweqwq}. \hfill $\blacksquare$

\section{Proof for non-convergence of generators}\label{a3}
{\bf Proof of Theorem \ref{th2001}:}
We have
\begin{multline}\label{eq301}
\parallel (H - H_\epsilon)\langle f, \cdot \rangle
\parallel^2_{L^2(\mu_\epsilon)} =
\int_{\Gamma_\epsilon}H \langle f, \omega \rangle H
\langle f, \omega \rangle \,d\mu_\epsilon(\omega) \\
- 2 \int_{\Gamma_\epsilon} H \langle f, \omega \rangle H_\epsilon
\langle f, \omega \rangle \,d\mu_\epsilon(\omega)
+ \int_{\Gamma_\epsilon}H_\epsilon \langle f, \omega \rangle H_\epsilon
\langle f, \omega \rangle \,d\mu_\epsilon(\omega).
\end{multline}
A direct consequence of Theorem \ref{th9}{\rm (ii)} is that
\begin{eqnarray*}
\lim_{\epsilon \to 0} \int_{\Gamma_\epsilon}H \langle f, \omega \rangle H
\langle f, \omega \rangle \,d\mu_\epsilon(\omega) =
\frac{(\rho_{\phi}^{(1)}(\beta, 1))^2}{\chi_{\phi}(\beta) } \parallel \Delta f
\parallel^2_{L^2(dx)}.
\end{eqnarray*}
Furthermore, we have
\begin{eqnarray}\label{eq355}
\lim_{\epsilon \to 0} \int_{\Gamma_\epsilon}H \langle f, \omega
\rangle H_\epsilon \langle f, \omega \rangle
\,d\mu_\epsilon(\omega) = \frac{(\rho_{\phi}^{(1)}(\beta,
1))^2}{\chi_{\phi}(\beta) } \parallel \Delta f
\parallel^2_{L^2(dx)},
\end{eqnarray}
where (\ref{eq355}) can be shown in the same way
as the convergence of Dirichlet forms in Theorem \ref{th11},
the argument is even more simple.
Showing convergence of the third term in (\ref{eq301}) is a quite
elaborate task. Using (\ref{eq25}), the coercivity identity
provided in Theorem \ref{dfdrr}, (\ref{eq302}) and Lemma \ref{le7} we obtain
\begin{multline*}
\int_{\Gamma_\epsilon}H_\epsilon \langle f, \omega \rangle
H_\epsilon \langle f, \omega \rangle \,d\mu_\epsilon(\omega) =
\epsilon^d \int_{\Gamma} H_{\tilde{\mu}_\epsilon} \langle f,
\gamma \rangle H_{\tilde{\mu}_\epsilon} \langle f, \gamma \rangle
\,d\tilde{\mu}_\epsilon(\gamma) \\ =
\rho_{\phi_\epsilon}^{(1)}(\beta,\epsilon^{-d}) \parallel \Delta f
\parallel^2_{L^2(dx)} + \, \frac{\epsilon^{-d}}{2} \, \beta \int_{{\mathbb
R}^{2d}} \Big{(} \nabla f(x) - \nabla f(y), \\
\nabla^2\phi_{\epsilon}(x-y)(\nabla f(x) - \nabla f(y))
\Big{)}_{{\mathbb R}^d} \,\,
\rho_{\phi_\epsilon}^{(2)}(\beta,\epsilon^{-d}, x, y) \,dx \,dy  \\
= \rho_{\phi}^{(1)}(\beta,1) \parallel \Delta f
\parallel^2_{L^2(dx)} + \, \frac{\epsilon^{-(d+2)}}{2} \, \beta
\int_{{\mathbb R}^{2d}} \Big{(} \nabla f(x) - \nabla f(y), \\
\nabla^2\phi\Big{(} \frac{x-y}{\epsilon}\Big{)} ( \nabla f(x) -
\nabla f(y)) \Big{)}_{{\mathbb R}^d} \,\, \rho_{\phi}^{(2)}\Big{(}
\beta, 1, \frac{x-y}{\epsilon}, 0 \Big{)} \,dx \,dy.
\end{multline*}
By the mean value theorem, we get
\begin{multline*}
\frac{\epsilon^{-(d+2)}}{2} \Big{(} \nabla f(x) - \nabla f(y),
\nabla^2\phi\Big{(} \frac{x-y}{\epsilon}\Big{)}(\nabla f(x)
- \nabla f(y)) \Big{)}_{{\mathbb R}^d} \\ =
\frac{\epsilon^{-d}}{2} \int_0^1 \int_0^1 \Big{(} \nabla^2
f(y + q_1(x-y)) \frac{x-y}{\epsilon}, \nabla^2\phi\Big{(}
\frac{x-y}{\epsilon}\Big{)} \nabla^2f(y + q_2(x-y))
\frac{x-y}{\epsilon} \Big{)}_{{\mathbb R}^d} \,dq_1 \,dq_2
\end{multline*}
Thus, we obtain an approximate identity and
\begin{multline*}
\lim_{\epsilon \to 0} \int_{\Gamma_\epsilon}H_\epsilon \langle
f, \omega \rangle H_\epsilon \langle f, \omega \rangle
\,d\mu_\epsilon(\omega) =  \rho_{\phi}^{(1)}(\beta,1) \parallel \Delta f
\parallel^2_{L^2(dx)} \\
+ \frac{1}{2} \sum_{i, j, k, l = 1}^d
\int_{{\mathbb R}^{d}} \beta \, x^k
x^l \, \partial_i \partial_j \phi(x) \,
\rho_{\phi}^{(2)}(\beta, 1, x, 0) \,dx \, \int_{{\mathbb
R}^{d}} \partial_i \partial_k f(y) \, \partial_j \partial_l f(y) \,dy,
\end{multline*}
where $x^i$ is the $i$-th component of $x \in {\mathbb R}^d$. Set
\begin{eqnarray}\label{eq573}
D_\phi(\beta, i, j, k, l) := \int_{{\mathbb R}^{d}} \beta
\, x^k x^l \, \partial_i \partial_j \phi(x)
\, \rho_{\phi}^{(2)}(\beta, 1, x, 0) \,dx.
\end{eqnarray}
For isotropic potentials the corresponding second correlation
function is also isotropic and the coefficient (\ref{eq573}) turns out to
be:
\begin{eqnarray*}
\int_{{\mathbb R}^{d}}
\beta \Big{(} \frac{x^{i} x^{j}}{r^3} (r V''(r) - V'(r)) +
\frac{V'(r)}{r} \delta_{i,j}\Big{)} \, x^k x^l
\, \tilde{\rho}_{\phi}^{(2)}(\beta, 1, r) \,dx,
\end{eqnarray*}
here $\delta_{i,j}$ is the Kronecker delta. Hence, for isotropic
potentials the coefficient $D_\phi(\beta, i, j, k, l)$ is only different from
zero if each index in the set $\{i,j,k,l \}$ at least occurs
twice. Utilizing polar coordinates, the identity
$\int_0^{2 \pi} \sin^4(\theta) \,d\theta = 3 \int_0^{2 \pi}
\sin^2(\theta) \cos^2({\theta}) \,d\theta$
and the symmetry of
$\int_{{\mathbb R}^{d}} \partial_i \partial_k f(y) \, \partial_j \partial_l f(y)
\,dy$
in all its indexes, we find
\begin{eqnarray*}
&& \lim_{\epsilon \to 0} \int_{\Gamma_\epsilon}H_\epsilon \langle
f, \omega \rangle H_\epsilon \langle f, \omega \rangle
\,d\mu_\epsilon(\omega) =  D_\phi(\beta) \parallel \Delta f
\parallel^2_{L^2(dx)}
\end{eqnarray*}
where
\begin{eqnarray*}
D_\phi(\beta) = \rho_{\phi}^{(1)}(\beta,1) + \frac{1}{2}
\int_{{\mathbb R}^{d}} \beta \, x^1 x^1 \, \partial_1 \partial_1 \phi(x) \,
\rho_{\phi}^{(2)}(\beta, 1, x, 0) \,dx.
\end{eqnarray*}

Next, we compare the coefficients $D_\phi(\beta)$ and
$(\rho^{(1)}_{\phi}(\beta, 1))^2/{\chi_{\phi}(\beta)}$ in terms of
a high temperature expansion. The latter coefficient is the
{\it isothermal compressibility} of the fluid or gas characterized by
$\mu$. Applying Theorem \ref{pr200}, we obtain
\begin{align*}
D_\phi (\beta) & =  1 + \beta^2 \frac{\partial^2}{\partial
\beta^2}D_\phi(0)+ o(\beta^2),  \\ \frac{(\rho^{(1)}_{\phi}(\beta,
1))^2}{\chi_{\phi}(\beta)} & =  1 + \beta^2
\frac{\partial^2}{\partial \beta^2}
\frac{(\rho^{(1)}_{\phi}(\cdot, 1))^2}{\chi_{\phi}}(0) +
o(\beta^2), \qquad \beta \in [0, \beta_0],
\end{align*}
where
\begin{align*}
\frac{\partial^2}{\partial \beta^2}D_\phi(0) & =
- \Big(\int_{{\mathbb R}^d} \phi(x) \,dx \Big)^2 +
\int_{{\mathbb R}^d} (x^1 \partial_1 \phi(x))^2 \,dx, \\
\frac{\partial^2}{\partial \beta^2}
\frac{(\rho^{(1)}_{\phi}(\cdot, 1))^2}{\chi_{\phi}}(0) & =
- \Big(\int_{{\mathbb R}^d} \phi(x) \,dx \Big)^2.
\end{align*}
Thus, the remainder function is given by
\begin{multline}\label{eq0524}
R_\phi(\beta) = D_\phi(\beta) -
\frac{(\rho^{(1)}_{\phi}(\beta, 1))^2}{\chi_{\phi}(\beta)}
= \beta^2 \int_{{\mathbb R}^d} (x^1 \partial_1 \phi(x))^2 \,dx
+ o(\beta^2), \, \beta \in [0, \beta_0].
\end{multline}
Hence, $R_\phi \equiv 0$ on $[0, \beta_0]$ is equivalent to
$\int_{{\mathbb R}^d} (x^1 \partial_1 \phi(x))^2 \,dx = 0$.
This, in turn, is equivalent to
\begin{gather}\label{eq83}
\partial_1 \phi(x) = 0 \quad  \mbox{for} \quad dx\mbox{-a.e.}
\quad x \in {\mathbb R}^d.
\end{gather}
Since the potential is isotropic, the only potential in
consideration which fulfills (\ref{eq83}) is $\phi \equiv 0$.
Hence, by (\ref{eq0524}) for $\mu \neq \pi_1$ there exist $\beta_1
\in (0, \beta_0]$ such that $R_\phi(\beta) > 0$
for all $\beta \in (0, \beta_1]$.
\hfill $\blacksquare$
\end{appendix}

\end{document}